\documentclass[11pt]{article}
\usepackage{amsfonts}
\usepackage{latexsym}

\newcommand{\qed}{$\hfill\Box$}

\newcommand{\R}{\mathbb{R}}
\newcommand{\LL}{\mathbb{L}}
\newcommand{\E}{\mathbb{E}}
\newcommand{\PP}{\mathbb{P}}
\newcommand{\N}{\mathbb{N}}

\newcommand{\al}{\alpha}
\newcommand{\be}{\beta}
\newcommand{\la}{\lambda}
\newcommand{\ga}{\gamma}

\newcommand{\Ga}{\Gamma}
\newcommand{\si}{\sigma}
\newcommand{\Si}{\Sigma}
\newcommand{\ep}{\varepsilon}

\newcommand{\De}{\Delta}
\newcommand{\om}{\omega}
\newcommand{\Om}{\Omega}
\newcommand{\ze}{\zeta}
\newcommand{\vp}{\varphi}

\newcommand{\rn}{\sqrt{n}}

\def\n{\hbox{$\mathcal N$}}

\def\ea{\hbox{$\mathcal E$}}
\def\g{\hbox{$\mathcal G$}}

\def\f{\hbox{$\mathcal F$}}

\def\ll{\hbox{$\mathcal L$}}

\def\nib{\noindent\bf}
\def\ni{\noindent}

\newcommand{\dd}{\Delta^n_i}

\def\is{\hbox{$\star (\mu -\nu)$}}

\def\proba{\hbox{$(\Omega ,\f,(\f_t)_{t\geq0},\PP)$}}
\def\prob'{\hbox{$(\Omega ',\f ',(\f '_t)_{t\geq0})$}}

\def\top{\hbox{$~\stackrel{\PP}{\longrightarrow}~$}}

\def\tol{\hbox{$~\stackrel{\ll}{\longrightarrow}~$}}
\def\tosc{\hbox{$~\stackrel{\mbox{Sk}}{\longrightarrow}~$}}
\def\toucp{\hbox{$~\stackrel{\mbox{\tiny u.c.p.}}{\longrightarrow}~$}}

\def\BV{\overline{V}}

\def\BZ{\overline{Z}}

\def\Bb{\overline{b}}

\newcommand{\PV}{\Pi}

\newcommand{\vsc}{\vskip 5mm}

\newcommand{\vst}{\vskip 3mm}
\newcommand{\vsq}{\vskip 4mm}

\newcommand{\bee}{\begin{equation}}
\newcommand{\eee}{\end{equation}}
\newcommand{\bea}{\begin{eqnarray}}
\newcommand{\eea}{\end{eqnarray}}
\newcommand{\bean}{\begin{eqnarray*}}
\newcommand{\eean}{\end{eqnarray*}}

\renewcommand{\theequation}{\arabic{section}.\arabic{equation}}

\newtheorem{prop}{Proposition}[section]

\newtheorem{lem}[prop]{Lemma}
\newtheorem{theo}[prop]{Theorem}

\begin{document}

\title{Asymptotic properties of power variations of L\'evy processes}
\author{Jean Jacod
\thanks{
Institut de Math\'ematiques de Jussieu, 175 rue du Chevaleret
75 013 Paris, France (CNRS -- UMR 7586,
and Universit\'e Pierre et Marie Curie--P6)}}

\date{}

\maketitle

\vst

\vst \noindent \it Institut de Math\'ematiques de Jussieu, \\
175 rue du
Chevaleret, 75 013 Paris, France \\
(CNRS -- UMR 7586,
and Universit\'e Pierre et Marie Curie--P6)\\
 e--mail:
jj@ccr.jussieu.fr \rm \vskip2cm

\rm\small

We determine the asymptotic behavior of the realized power
variations, and more generally of sums of a given function $f$
evaluated at the successive increments of a L\'evy process between
the successive times $i\De_n$ for $i=0,1,\ldots,n$. One can elucidate
completely the first order behavior, that is the convergence in
probability of such sums, possibly after normalization and/or
centering: it turns out that there is a rather wide variety of
possible behaviors, depending on the structure of jumps and on
the chosen test function $f$. As for the associated central limit theorem,
one can show some versions of it, but unfortunately in a limited
number of cases only: in some other cases a CLT just does not
exist.

\normalsize

\vst


\newpage

\section{Introduction}\label{sec-Intro}
\setcounter{equation}{0}
\renewcommand{\theequation}{\thesection.\arabic{equation}}

The behavior of the so--called realized power or multipower
variations for discretely observed semimartingales has attracted a lot
of attention
recently, mainly because of applications in finance: they are used
for estimating the volatility or integrated volatility, typically in a
stochastic volatility context, with or without jumps. They are useful
also when one wants to
decide whether a discretely observed process is continuous or has
jumps. In all these cases one is looking at what kind of power or
multipower or possibly truncated power variations is best for
estimating the integrated volatility and/or some functions of the jumps.
Let us
quote for example \cite{ABD}, \cite{BS}, \cite{BGJPS}, and some
partial results when the process under considerations has jumps
may be found in \cite{BSW}, \cite{Ma} or \cite{W}, see also
\cite{AJ} in a different context.

Determining the asymptotic behavior of power or multipower variations
when the time lag goes to $0$ is not a trivial problem,
especially when there are jumps. In this paper we consider the
simplest possible case: the underlying process $X$ is a
L\'evy process, and we look at power or truncated power variations
only. However, one may hope that the results extend to more general
semimartingales, at least for those
having ``absolutely continuous'' characteristics. They also
extend to multipower variations, but probably not in a trivial way.

When $X$ is a continuous L\'evy process the results are of course elementary,
but they become not so simple and indeed quite
versatile when there are jumps. However the homogeneity and the
independent increments property of $X$ allow to give complete answers
for the convergence in probability of (possibly normalized) power
variations to a non trivial limit. As for the associated CLT we give
a complete answer only when the Gaussian part of $X$ does not vanish,
and this probably the interesting case in practice. When the Gaussian
part vanishes, the existence of a CLT depends on the behavior of the
L\'evy measure of $X$ near $0$, and such a CLT does not necessarily
exist: in this paper we essentially do not consider this situation.

Let us be more specific. We have a $1$--dimensional L\'evy process
$X$ with characteristics $(b,c,F)$ (see Section \ref{sec2} for
more details). This process starts at $0$  (i.e. $X_0=0$) and it
is sampled at times $\De_n,2\De_n,\cdots$, where $\De_n$ is {\em
always} assumed to go to $0$ as $n\to\infty$. With the notation
\bee\label{1} \dd Y=Y_{i\De_n}-Y_{(i-1)\De_n} \eee for the
increments of any process $Y$, we consider the {\em realized
$r$--variation process} \bee\label{2}
\PV^n(r)_t:=\sum_{i=1}^{[t/\De_n]}|\dd X|^r \eee where
$r\in(0,\infty)$, and also the truncated realized $r$--variation
at level $a$: \bee\label{3}
\PV^n(r,a)_t:=\sum_{i=1}^{[t/\De_n]}|\dd X|^r1_{\{|\De^n_iX|\leq
a\}}. \eee We are interested in the asymptotic behavior of
$\PV^n(r)$, or of $\PV^n(r,a_n)$ for a sequence of truncation
levels $a_n$ of the form $a_n=a$ or $a_n=a\sqrt{\De_n}$ for some
$a>0$, although other choices for $a_n$ would also be possible;
the choice $a_n=a\sqrt{\De_n}$ is obviously related to the
property that $\dd X/\sqrt{\De_n}$ is $\n(0,c)$ when $X$ is
Gaussian, and it reflects the fact that we are mainly interested
in the case where $c>0$. The LLN results (that is, the convergence
in probability, possibly after normalization) for $\PV^n(r)$ cover
all possible situations, but for the truncated version
$\PV^n(r,a)$ with a fixed level $a$ there are difficulties due to
the non--continuity of the function $x\mapsto |x|^r1_{\{|x|\leq
a\}}$.

From a statistical point of view, it is also natural to consider that
we sample the process $X$
over the interval $[0,T_n]$, where $T_n=n\De_n$. Then we are
interested in the behavior of the variables $\PV^n(r)_{T_n}$, or more
generally of the processes
\bee\label{4}
\overline{\PV}^n(r)_t=\sum_{i=1}^{[nt]}|\dd X|^r,\qquad
\overline{\PV}^n(r)_t=\sum_{i=1}^{[nt]}|\dd X|^r1_{\{|\De^n_iX|\leq a\}},
\eee
so $\PV^n(r,a)_{T_n}=\overline{\PV}^n(r,a)_1$. In practice two cases are
of interest: one is when $T_n=T$ does not depend on $n$ (that is
$\De_n=T/n$ for some $T>0$), and in this case $\PV^n(r,a)$ and
$\overline{\PV}^n(r,a)$ are essentially identical since
$\overline{\PV}^n(r,a)_t=\PV^n(r,a)_{t/T}$. The other case is when
$T_n\to\infty$: the results are then much simpler than for $\PV^n(r,a)$,
but they require some integrability which holds for the truncated
version, but not in general for $\overline{\PV}^n(r)_t$.
\vst

The two afore--mentioned difficulties -- lack of continuity and
lack of integrability -- lead us to consider slightly
more general processes, not really any more difficult to study. More
precisely, for any function $f$ on $\R$ we set
\bee\label{5}
V^n(f)_t:=\sum_{i=1}^{[t/\De_n]}f(\dd X),\qquad
\BV^n(f)_t:=\sum_{i=1}^{[nt]}f(\dd X).
\eee
Then $\PV^n(r,a)=V^n(f)$ with $f(x)=|x|^r1_{\{|x|\leq a\}}$ (for
$a\in(0,\infty]$), and accordingly for $\overline{\PV}^n(r,a)$. The
truncated case with a truncation level $a\sqrt{\De_n}$ is
accommodated by looking at the following modified processes:
\bee\label{6}
V'^n(f)_t:=\sum_{i=1}^{[t/\De_n]}f(\dd X/\sqrt{\De_n}),\qquad
\BV'^n(f)_t:=\sum_{i=1}^{[nt]}f(\dd X/\sqrt{\De_n}).
\eee
So $\PV^n(r,a\sqrt{\De_n})=\De_n^{r/2}V'^n(f)$, with
$f(x)=|x|^r1_{\{|x|\leq a\}}$ again.

For the processes $V^n(f)$ and $\BV^n(f)$, it turns out that the
crucial property of $f$ is
its behavior near $0$. Depending on the case at hand, we will restrict
our attention to functions which, for some $r>0$, are ~o$(|x|^r)$ or
~O$(|x|^r)$ as $x\to0$, or even which coincide
with the function $x\mapsto|x|^r$ on a neighborhood of $0$.

For the processes $V'^n(f)$ and $\BV'^n(f)$, however, the whole
function $f$ is important and the behavior of $f$ near $0$ does
not play a specific role. On the other hand the special truncation
we consider (at $\sqrt{\De_n}$) leads us to give result only when
the Gaussian part of $X$ is not vanishing. If it is, other
''truncation rates'' would be more appropriate.

The paper is organized as follows: in Section 2 we set the notation
and state the results
pertaining to the processes $V^n(f)$ and $V'^n(f)$. Our concern of
being as complete as possible leads to a rather large number of
results corresponding to various situations, which are all encountered
for the processes $\PV^n(r,a)$ already.
Section 3 is devoted to stating the results for the (less
important) processes $\BV^n(f)$ and $\BV'^n(f)$.
The proofs are given in Sections 4--7.

\section{The processes $V^n(f)$ and $V'^n(f)$}\label{sec2}
\setcounter{equation}{0}
\renewcommand{\theequation}{\thesection.\arabic{equation}}

\subsection{General setting and notation.}\label{SN}

Our L\'evy process $X$ is defined on the filtered probability space
$\proba$ and we assume that $(\f_t)$ is the right--continuous
filtration generated by $X$ and that $\f=\bigvee\f_t$. It has characteristics
$(b,c,F)$: this is the unique triple $b\in\R$, and $c\geq0$, and $F$ a
measure on $\R$ which does not charge $\{0\}$ and integrates $x\mapsto
x^2\bigwedge1$, such that
\bee\label{Levy}
\E(e^{iuX_t})=\exp t\left(iub-\frac{cu^2}2+
\int F(dx)\left(e^{iux}-1-iux1_{\{|x|\leq1\}}\right)\right).
\eee
We denote by $X^c$ the continuous martingale part (the ``Gaussian
part'') of $X$, and also by $\mu$ the jump measure of $X$ and by $\nu$
the predictable compensator of $\mu$, that is $\nu(ds,dx)=dsF(dx)$.
We use the symbol $\star$ to denote the stochastic integral with
respect to the measures $\mu$ and $\mu-\nu$ (cf. \cite{JS} for more
details on this notation and the decomposition below). Then
\bee\label{10}
X_t=bt+X^c_t+(x1_{\{|x|\leq1\}})\star(\mu-\nu)_t+(x1_{\{|x|>1\}})\star\mu_t.
\eee
Recall that $X^c=\si~W$, where $\si=\sqrt{c}$ and $W$ is standard
Brownian motion.

For any function $f$ for which the following makes sense, we write
\bee\label{9}
F(f)=\int f(x)F(dx),\qquad  H_t(f)=\E(f(X_t)),\qquad
\Ga_t(f)=H_t(f^2)-H_t(f)^2.
\eee
The following functions $\phi_r$ will often occur:
\bee\label{99}
\phi_r(x)~=~\left\{\begin{array}{ll}
1\bigwedge|x|^r\qquad&\mbox{if }~~0<r<\infty\\
1 &\mbox{if }~~r=0,\end{array}\right.
\eee
and the set
\bee\label{13}
I=\{r\geq0:F(\phi_r)<\infty\}
\eee
plays a fundamental role. It is an interval of the form
$[\al,\infty)$ or $(\al,\infty)$, for some $\al\in[0,2]$. We
have $2\in I$ always, end the process $X-X^c$ is locally of finite
variation if and only if $1\in I$, and $X$ has locally finitely many
jumps if and only if $0\in I$. Set
\bee\label{14}
\left.\begin{array}{l} \Bb=b-\int_{\{|x|\leq1\}}xF(dx),\quad
X'_t=\sum_{s\leq t}\De X_s
\qquad\mbox{if }~1\in I\\[2mm]
X''=X-X^c.\end{array}\right\} \eee So if $1\in I$ we have
$X''_t=\Bb t+X'_t$, and $\Bb$ is the ``genuine'' drift.

Let us introduce several classes of functions on $\R$. First, for
$r\in(0,\infty)$ we denote by $\ea_r$ and $\ea'_r$
and $\ea''_r$ the following sets of Borel functions;
\bee\label{69}
\left.\begin{array}{lll}
\ea_r&:&\mbox{all $f$ with polynomial growth and
~$f(x)\sim|x|^r$~ as $x\to0$}\\
\ea'_r&:&\mbox{all $f$ locally bounded and with ~$f(x)=$ O$(|x|^r)$~
as $x\to0$}\\
\ea''_r&:&\mbox{all $f$ locally bounded and with ~$f(x)=$ o$(|x|^r)$~
as $x\to0$}.\end{array}\right\}
\eee
We write $\ea_r^b$, $\ea'^b_r$ and $\ea''^b_r$ for the sets of bounded
functions belonging to $\ea_r$, $\ea'_r$ and $\ea''_r$ respectively.
Finally, as usual  $C^p$ denotes
the set of $p$ times continuously differentiable functions (for
$p\geq0$), and $C^{0,F}$ is the set of functions that are $F$--a.e.
continuous. We have $\ea_r\subset\ea'_r\subset \ea''_s\subset\ea'_s$
for $s<r$, and $\phi_r\in\ea^b_r\cap C^0$.

If $f\in\ea'_r$ then $(|f|\bigwedge1)*\nu$ and thus $f*\mu$ as well are
finite--valued as soon as $r\in I$, whereas if $r\notin I$ then
$f*\mu_t=\infty$ a.s. for all $t>0$ and $f\in\ea_r$ (to see this, observe
that if further $f\geq0$, the Laplace transform of $f\star\mu_t$ is
$\la\mapsto \exp-t
\int_{\{|x|\leq1\}}\left(1-e^{-\la f(x)}\right)F(dx)=0$ for all $\la$).
In connection with this one can also define the process
\bee\label{15}
\Si(f,\vp)=(f\vp)\is+(f(1-\vp))\star\mu.
\eee
Here, $\vp$ denotes a function having $1_{[-a,a]}\leq\vp\leq1$ for some
$a>0$ and having a compact support when $f$ is unbounded; the last
integral above is always meaningful, and the
first (stochastic) integral makes sense as soon as $f\in\ea'_r$ for some $r$
such that $2r\in I$.

Finally we give some miscellaneous notation. We denote by $U$ a
generic $\n(0,1)$ variable, and by
$\mu_r$ the $r$th absolute moment of $U$. With any process $Y$ we associate
the increments $\dd Y$ by (\ref{1}) and also the ''discretized process''
\bee\label{1'}
Y^{(n)}_t~=~Y_{\De_n[t/\De_n]}~=~\sum_{i=1}^{[t/\De_n]}\dd Y.
\eee
We say
that a sequence of processes $Y^n$ ''converges u.c.p.'' to a process $Y$,
denoted $Y^n\toucp Y$ or $Y^n_t\toucp Y_t$, if $\sup_{s\leq
t}|Y_s^n-Y_s|\top 0$ for all $t>0$.

Recall also (see e.g. \cite{JS})
that the stable convergence in law of $Y^n$ to $Y$, when each $Y^n$ is
a c\`adl\`ag process on $(\Om,\f,\PP)$ and $\f$ is the $\si$--field
generated by all $X_t$'s and $Y$ is a c\`adl\`ag process defined on an
extension of the space $(\Om,\f,\PP)$, means that
$\E(Zg(Y^n))\to\E(Zg(Y))$ for all bounded $\f$--measurable variable
$Z$ and all bounded continuous function $g$ on the space of all
c\`adl\`ag functions, endowed with the Skorokhod topology.

\subsection{Results: the law of large numbers.}\label{res}

\begin{theo}\label{T3} (i) The processes $V^n(f)$   converge in
  probability to a limit $V(f)$, for the Skorokhod topology, in the
  following cases:

\noindent(a) With $V(f)=f\star\mu$, when

\hskip2cm [a-1] $f\in\ea_2''\cap C^{0,F}$,

\hskip2cm[a-2] $f\in\ea_r'\cap C^{0,F}$ if $r\in I\cap(0,1)$ and $c=0$,

\hskip2cm[a-3] $f\in\ea_1''\cap C^{0,F}$ if $1\in I$ and $c=0$,

\hskip2cm[a-4] $f\in\ea_r'\cap C^{0,F}$ if $r\in I\cap(0,1]$ and
$c=\Bb=0$.

\noindent(b) With $V(f)_t=f\star\mu_t+ct$, when $f\in\ea_2\cap
C^{0,F}$.

\noindent(c) With $V(f)_t=f\star\mu_t+|\Bb|t$, when $f\in\ea_1\cap C^{0,F}$
and $c=0$ and $1\in I$.

(ii) In particular this holds for $f\in\ea_r\cap C^{0,F}$ when
$r\in[2,\infty)$, or
when $r\in I\cap[1,2)$ if $c=0$, or when $r\in I\cap(0,1)$ if
$c=0$ and $\Bb=0$.

(iii) If $r$ is not as in (ii), then
$\frac1{\De_n}~H_{\De_n}(\phi_r)\to\infty$ and, for any $f\in\ea_r$, the
processes $\frac{\De_n}{H_{\De_n}(\phi_r)}~V^n(f)_t$ converge
u.c.p. to $t$; in particular, $V^n(f)_t\top+\infty$ for all $t>0$.
\end{theo}

When $f(x)=|x|^r$ the case (b) ($r=2$) is of course well known (this is the
convergence of the realized quadratic variation), and (a) for $r>2$ may be
found in \cite{L} for general semimartingales, and (c) ($r=1$) is also well
known because then $V(f)$ is the variation process of $X$.

(iii) is not really satisfactory, since the rate
$\frac1{\De_n}~H_{\De_n}(\phi_r)$ is not ``explicit''. In the most
interesting case $c>0$ it is however possible to give a more explicit
result. In this case we can also provide an LLN for $V'^n(f)$~:

\begin{theo}\label{T4} Assume that $c>0$.

(i) If $f\in\ea_r$ for some $r\in(0,2)$ then
$\De_n^{1-r/2}V^n(f)_t$ converges u.c.p. to $c^{r/2}\mu_rt$. The same
holds for $r\in[2,\infty)$ if there is no jump (i.e. $F$=0).

(ii) If $f$ is bounded then
$\De_nV'^n(f)_t$ converges u.c.p. to $t\E(f(\si U))$.
\end{theo}

Note that, apart from the condition $r<2$ in (i) when $F\neq 0$, these
results do not depend on the jumps at all, and moreover in (i) the
limit does not even depend on the function $f$ (subject to being in
$\ea_r$). Clearly (i) fails if we assume only $f\in\ea'_r$.

We also have another LLN for $V^n(f)$ without normalization in some of
the cases covered by (iii) of Theorem \ref{T3}, provided we center
$V^n(f)$ appropriately:

\begin{theo}\label{T5}  Let either $1<r<2$, or $r\leq1$ and $2r\in I$ and
  $c=0$. If~ $f\in\ea^b_r\cap C^{0,F}$ and $\vp\equiv1$, or if~
  $f\in\ea_r\cap C^{0,F}$ and $\vp$ is continuous with
$1_{[-a,a]}\leq\vp\leq1_{[-a',a']}$ for some $0<a<a'<\infty$, the processes
$V^n(f)_t-\frac1{\De_n}~H_{\De_n}(f\vp)t$ converge in probability, for
the Skorokhod topology, to $\Si(f,\vp)$.
\end{theo}

\ni\it Remark 1~: \rm This theorem overlaps with (i) of Theorem \ref{T3},
but in the
overlapping cases $\frac1{\De_n}~H_{\De_n}(f\vp)$ converges to
$F(f\vp)+|\Bb|$ when $r=1$ (and $c=0$), and to $F(f\vp)$
otherwise, so the two results are consistent. \qed
\vsq

\ni\it Remark 2~: \rm When $c>0$ and $1<r<2$ and $f\in\ea_r\cap
C^{0,F}$, we can apply
the two previous theorems at once: one can
rewrite (i) of Theorem \ref{T4} as
$\De_n^{1-r/2}\left(V^n(f)_t-\frac t{\De_n}~H_{\De_n}(f)\right)$
$\toucp 0$ (see Lemma \ref{L2} below), so one can view Theorem \ref{T5} in this
case as the CLT associated with the LLN of Theorem \ref{T4}, or
perhaps rather as a ''second order'' LLN because the convergence
takes place in probability. \qed
\vsq

\ni\it Remark 3~: \rm One could also prove that the result holds for
$f$ in $\ea'_r$ instead of $\ea_r$. This is in contrast with Theorem
\ref{T6}--(ii) below, which is analogous to the above in a sense, but
for which we cannot substitute $\ea_1$ with $\ea'_1$. \qed

\subsection{Results: the CLT.}\label{res-CLT}

We now turn to the CLT, for which various versions occur. First
there are CLTs associated with LLNs having a deterministic limit
(possibly after centering): for $V'^n(f)$ that concerns Theorem
\ref{T4}--(ii), and for $V^n(f)$ for $f\ea_r$ that concerns the cases of (iii)
of Theorem \ref{T3} not covered by Theorem \ref{T5}, that is

(1)~~ $c=0$, $r<1$ and $2r\notin I$;

(2)~~ $c>0$, $r=1$;

(3)~~ $c>0$, $r<1$.

\ni For (1) the answer is complete, but probably useless because
the rate is not explicit. For (2) the answer is half--way between
a genuine CLT and an LNN: these two cases are dealt with in
Theorem \ref{T6}. Case (3) is considered later, together with the
CLT about $V'^n(f)$.

Second, there are CLTs associated with LLN having random limits,
that is (i) of Theorem \ref{T3}. When $c>0$ we have such a CLT
below, which extends results of \cite{JP}; this concerns the
cases $r=2$ and $r>3$ only. When $c=0$ if a CLT exists it is very
complicated, and indeed it does not exist in general (see \cite{J}
for special examples of this situation): we do not consider those
cases here at all.

\begin{theo}\label{T6}
(i) Let $f\in\ea_r$ for some $r\in(0,1)$ and $c=0$ and $2r\notin I$.
Then $\frac1{\De_n}~\Ga_{\De_n}(\phi_r)\to\infty$, and the processes
$\frac{\De_n}{\sqrt{\Ga_{\De_n}(\phi_r)}}\left(V^{n}(f)_t-
\frac t{\De_n}H_{\De_n}(\phi_r)\right)$ converge stably in law to a
standard Brownian motion independent of $X$.

(ii) If $c>0$ and either $f\in\ea_1^b\cap C^{0,F}$ and $\vp=1$, or
$f\in\ea_1\cap C^{0,F}$ and $\vp$ is continuous with
$1_{[-a,a]}\leq\vp\leq1_{[-a',a']}$ for some $0<a<a'<\infty$, the processes
$V^n(f)_t-\frac t{\De_n}~H_{\De_n}(f\vp)$ converge stably in law to
the process $\Si(f,\vp)+W'$, where $W'$ is a Wiener process
independent of $X$ and with unit variance $c(\mu_2-\mu_1^2)
=c(1-2/\pi)$.
\end{theo}

\ni\it Remark 4~: \rm
In the situation of (i) above, $\frac{\De_n}{H_{\De_n}(\phi_r)}~
\left(V^{n}(f)_t-\frac t{\De_n}H_{\De_n}(\phi_r)\right)\toucp0$ by
Theorem \ref{T3}--(ii), hence we also have
 $H_{\De_n}(\phi_r)/\sqrt{\Ga_{\De_n}(\phi_r)}\to\infty$. \qed
\vst

Next we give the CLT associated with Theorem \ref{T4}, with $r<1$
in (i) (case (3) above). Observe the two different centerings
below.

\begin{theo}\label{T7} Assume that $c>0$.

(i) Let $f\in\ea_r$ for some $r$ in $(0,1)$ in general, and
some $r\in(0,\infty)$ if there is no jump.

\ni 1) The processes
$\frac1{\sqrt{\De_n}}\left(\De_n^{1-r/2}V^n(f)_t-
t\De_n^{-r/2}H_{\De_n}(\phi_r)\right)$ converge stably in law to a Wie\-ner
process independent of $X$ and with unit variance
$c^r(\mu_{2r}-\mu_r^2)$,
and if further $f$ is bounded we can substitute
$H_{\De_n}(\phi_{r})$ with $H_{\De_n}(f_j)$.

\ni 2) If $1\in I$ the same result holds for the processes
$\frac1{\sqrt{\De_n}}\left(\De_n^{1-r/2}V^n(f)_t-
tc^{r/2}\mu_r\right)$, whereas
$\frac1{\De_n^{1-s/2}}\left(\De_n^{1-r/2}V^n(f)_t-
tc^{r/2}\mu_r\right)\toucp0$ if $1\notin I$
and $s\in I$ for some $s\in(1,2)$.

(ii) Let $f$ be a bounded function.

\ni 1) The processes
$\frac1{\sqrt{\De_n}}\left(\De_nV'^n(f)_t-
t\E(f(X_{\De_n}/\sqrt{\De_n}))\right)$ converge stably in law to a Wiener
process independent of $X$ and with unit variance the variance of
$f(\si U)$.

\ni 2) Assume further that $f$ is an even function. If $1\in I$ the same
result holds for the processes
$\frac1{\sqrt{\De_n}}\left(\De_nV'^n(f)_t-t\E(f(\si U))\right)$,
whereas $\frac1{\De_n^{1-s/2}}\left(\De_nV'^n(f)_t-
t\E(f(\si U))\right)\toucp0$
if $1\notin I$ and $s\in I$ for some $s\in(1,2)$.
\end{theo}

Apart from the restriction $r<1$ in (i), we thus see
that the presence of jumps does not affect the CLT. But it is
important to note that the centering is {\em not}\, the one
expected from Theorem \ref{T4} if there are jumps, unless $1\in I$.
The last claims in (2) of both (i) and (ii) are clearly not sharp,
but we do not know the exact
rate of convergence in those cases.

Finally we give a CLT associated with LLNs having random limits,
when $c>0$. For this we need some more notation. We consider a
sequence $(U_n)_{n\geq1}$ of i.i.d. $\n(0,1)$ variables and also a
standard Brownian motion $W'$ independent of that sequence, all of
these defined on an extension of $\proba$ and independent of $X$.
Denote also an enumeration of the jump times of $X$ by
$(T_n)_{n\geq1}$, where the $T_n$'s are stopping times (not
necessarily increasing in $n$, of course).

\begin{theo}\label{T8} a) Assume that $c>0$ and let $g\in\ea_1'$. Then
  the process
\bee\label{90}
Z(g)_t=\sum_{p:~T_p\leq t} \si g(\De X_{T_p})U_p.
\eee
is well defined, and when $F(g^2)<\infty$ it is a locally
square--integrable martingale (on the extended space)
w.r.t. the filtration generated by the triple $(X,Z,W')$.

b) If $f$ is $C^1$ over $\R$ and $C^2$ on a neighborhood of $0$, with
$f''(x)=$ o$(|x|)$ as $x\to0$, the processes
$\frac1{\sqrt{\De_n}}~(V^n(f)-V(f)^{(n)})$ converge stably in law
toward $Z(f')$.

c) If $f\in C^1$ and $f(x)=x^2$ on a neighborhood of $0$, then the processes
$\frac1{\sqrt{\De_n}}~(V^n(f)-V(f)^{(n)})$ converge stably in law toward
$Z(f')+c\sqrt{2}~W'$.
\end{theo}

\ni\it Remark 5~: \rm
We do not have convergence in law of $\frac1{\sqrt{\De_n}}~
(V^n(f)-V(f))$, because of a peculiarity of Skorokhod topology:
indeed if $X$ has a jump at some time $S$, then the above process
has a ''big'' jump at $S$, and another ''big'' one at
$\De_n[S/\De_n]$, and those two times are close but distinct. But
if we substitute $V(f)$ with $V(f)^{(n)}$
then those two jumps are put together at the single time
$\De_n[S/\De_n]$, thus Skorokhod convergence may, and indeed
does, occur. However the processes $\frac1{\sqrt{\De_n}}~
(V^n(f)-V(f))$ converge finite--dimensionally, stably in law, to
the limit given above. \qed
\vsq

\ni\it Remark 6~: \rm In (b) the assumption is a bit more than asking
$f\in\ea''_3\cap C^1$. This theorem says
nothing when $f\in\ea_r\cap C^1$ for some $2<r\leq3$ and in these
cases the sequence $\frac1{\sqrt{\De_n}}~(V^n(f)-V(f)^{(n)})$ is tight
if $r=3$, but not tight if $2<r<3$.  To see this, suppose that
there is no jump and take $f\in\ea_r\cap C^{0,F}$; then the variables
$\frac1{\sqrt{\De_n}}~
(V^n(f)_t-V(f)^{(n)}_t)=\frac1{\sqrt{\De_n}}~V^n(f)_t$ converge in
probability to $+\infty$ when $2<r<3$, and to $c^{3/2}\mu_3t$ when
$r=3$, and to $0$ when $r>3$, by Theorem \ref{T4}. When
$r=2$, though, we have $V(f)_t=ct$ and then the previous theorem
reduces to Theorem \ref{T7}--(a), when $F=0$ and $r=2$.

One can also point out that when $f\in\ea_2\cap C^1$ but $f$ is not
equal to $|x|^r$ near $0$, we do not know whether the claim (c) holds true.
\qed
\vsq

Finally, for practical purposes we need a multidimensional CLT in
which we consider components as in (i) and (ii) of Theorem
\ref{T7} and also as in Theorem \ref{T8}. Below we have a
$d$--dimensional process and the index set $\{1,\ldots,d\}$ for
the components is partitioned into three (possibly empty) subsets
$J$, $J'$ and $J''$.

\begin{theo}\label{T9} Assume that $c>0$. Let $Y^n=(Y^{n,j})_{1\leq j\leq
d}$ be the process having the following components:

\ni $\bullet$ If $j\in J$, then $Y^{n,j}=
\frac1{\sqrt{\De_n}}\left(V^n(f_j)-V(f_j)^{(n)}\right)$, where
$f_j$ is like in Theorem \ref{T8}--(b),

\ni$\bullet$ If $j\in J'$, then $Y^{n,j}_t=
\frac1{\sqrt{\De_n}}\left(\De_n^{1-r(j)/2}V^n(f_j)_t-
t\De_n^{-r(j)/2}H_{\De_n}(\phi_{r(j)})\right)$, where $f_j\in\ea_{r(j)}$
for some $r(j)\in(0,1)$ in general or $r(j)\in(0,\infty)$ if there is
no jump,

\ni$\bullet$ If $j\in J''$, then $Y^{n,j}_t=
\frac1{\sqrt{\De_n}}\left(\De_nV'^n(f_j)_t-
t\E(f_j(X_{\De_n}/\sqrt{\De_n}))\right)$, where $f_j$ is bounded.

Then $Y^n$ converges stably in law to a process $Y$ whose components
$Y^j$ for $j\in J$ are of the form $Z(f'_j)$ given by (\ref{90})
with the same family $(U_p)$ for all $j$, and the components of
$Y$ in $J'\cup J''$ form a Wiener process independent of $X$ and
of the $(U_p)$'s and with unit variance--covariance given by
\bee\label{17}
c'^{jk}~=~\left\{\begin{array}{ll}
c^{r(j)/2+r(k)/2}(\mu_{r(j)+r(k)}-\mu_{r(j)}\mu_{r(k)}),
\quad&j,k\in J'\\[2mm]
\E(f_j(\si U)f_k(\si U))-\E(f_j(\si U))\E(f_k(\si U))
\quad&j,k\in J''\\[2mm]
\E(|\si U|^{r(j)}f_k(\si U))-c^{r(j)/2}\mu_{r(j)}\E(f_k(\si U))
\quad&j\in J',~k\in J''.
\end{array}\right.
\eee
Moreover we can replace $\De_n^{-r(j)/2}H_{\De_n}(\phi_{r(j)})$ by
$c^{r/2}\mu_r$ for $j\in J'$ and $\E(f_j(X_{\De_n}/\sqrt{\De_n}))$ by
$\E(f_j(\si U))$ for $j\in J''$, as soon as $1\in
I$ and provided also that $f_j$ is an even function for all $j\in J''$.
\end{theo}

\section{The processes $\BV^n(f)$ and $\BV'^n(f)$}\label{sec3}
\setcounter{equation}{0}
\renewcommand{\theequation}{\thesection.\arabic{equation}}

Now we give the results for the processes $\BV^n(f)$ and $\BV'^n(f)$
introduced in (\ref{5}) and (\ref{6}). As said in the introduction,
this is of interest only when $T_n=n\De_n$ goes to $\infty$, an
assumption made throughout.

The situation is much simpler than in the previous section, due
to the built--in ergodic properties of $X$, or rather of its
increments, and the proofs essentially boil down to using the usual LLN
and CLT: we never encounter random limits in the LLNs.
Moreover since $T_n\to\infty$ the stable
convergence in law is automatically implied by the mere convergence in
law, so we do not even state it in the results. The price we have to
pay is that the variable $f(X_{\De_n})$
should be integrable, or square--integrable. To avoid complicated
statements, we just suppose that $f$ is bounded. Also, because of this
necessary ergodicity, the results certainly cannot be extended to
general semimartingales, contrarily to what happens for the results of
the previous section.

We begin with the LLN, which corresponds to Theorems \ref{T3} and
\ref{T4}.

\begin{theo}\label{T10}  Assume $T_n\to\infty$.

(i) The processes $\frac1{T_n}~\!\!\BV^n(f)$ converge u.c.p. to a
limit $v(f)t$ in the following cases:

\noindent(a) With $v(f)=F(f)$, when $f$ is bounded and like in [a--1]--[a-4]
of Theorem \ref{Y3}

\noindent(b) With $v(f)=F(f)+c$, when $f\in\ea_2^b\cap
C^{0,F}$.

\noindent(c) With $v(f)=F(f)+|\Bb|$, when $f\in\ea_1^b\cap C^{0,F}$
and $c=0$ and $1\in I$.

(ii) If $c>0$ and $f\in\ea_r^b$ for some $r\in(0,2)$, the processes
$\frac1{n\De_n^{r/2}}~\BV^n(f)_t$ converge u.c.p. to
$c^{r/2}\mu_rt$. The same holds when $r\in[2,\infty)$ if there is no jump.

(iii) If $c>0$ and $f$ is bounded the processes $\frac1n~\BV'^n(f)$
converge u.c.p. to $t\E(f(\si U))$.
\end{theo}

As for the CLT, we have first a ``general'' result which covers all
situations but is hardly usable in practice because the rate in some
cases and the centering in most cases are not explicit:

\begin{theo}\label{T11} Assume $T_n\to\infty$.

(i) If either $f\in\ea''^b_1\cap C^{0,F}$, or $f\in\ea'_r\cap
C^{0,r}$ for some $r<1$ with $2r\in I$ and $c=0$, the processes
$\frac1{\sqrt{T_n}}\left(\BV^n(f)_t-nH_{\De_n}(f)t\right)$ converge
in law to a Wiener process with unit variance $F(f^2)$.

(ii) If $f\in\ea^b_1\cap C^{0,F}$ the processes
$\frac1{\sqrt{T_n}}\left(\BV^n(f)_t-nH_{\De_n}(f)t\right)$ converge
in law to a Wiener process with unit variance
$F(f^2)+c(\mu_2-\mu_1^2)=F(f^2)+c(1-2/\pi)$.

(iii) If $r<1$ and either $c>0$ or $c=0$ and $2r\notin I$, for all
$f\in\ea^b_r$ the processes
$\frac1{\sqrt{n\Ga_{\De_n}(f)}}\left(\BV^n(f)_t-nH_{\De_n}(f)t\right)$
converge in law to a standard Brownian motion.
\end{theo}

Finally we have a CLT similar to Theorem \ref{T7}, and for which we give the
multidimensional version only. Below the index set $\{1,\cdots,,d\}$
is partitioned into two subsets $J'$ and $J''$. Notice the difference
with Theorem \ref{9}: the first set $J$ of components is absent
here, although probably a result similar to Theorem \ref{T8} hold here
(but with a limit which is a Wiener process again, instead of a
discontinuous process).

\begin{theo}\label{T13} Assume $T_n\to\infty$ and $c>0$. Let
 $Y^n=(Y^{n,j})_{1\leq j\leq d}$ be the process having the following
 components:

\ni$\bullet$ If $j\in J'$, then $Y^{n,j}_t=
\rn\left(\frac1{n\De_n^{r(j)/2}}~\!\BV^n(f_j)_t-
t\De_n^{-r(j)/2}H_{\De_n}(\phi_{r(j)})\right)$, where $f_j\in\ea_{r(j)}$
for some $r(j)\in(0,1)$ in general or $r(j)\in(0,\infty)$ if there is
no jump,

\ni$\bullet$ If $j\in J''$, then $Y^{n,j}_t=
\rn\left(\frac1n~\!\BV'^n(f_j)_t-
t\E(f_j(X_{\De_n}/\sqrt{\De_n}))\right)$, where $f_j$ is bounded.

\ni Then $Y^n$ converges stably in law to a Wiener process with unit
variance--covariance given by (\ref{17}).

Moreover we can replace $\De_n^{-r(j)/2}H_{\De_n}(\phi_{r(j)})$ by
$c^{r/2}\mu_r$ for a $j\in J'$ it there is $s\in I\cap[0,1)$
such that the sequence $T_n\De_n^{(1-s)(1-r(j))}$ is bounded; in a
similar way, we can replace $\E(f_j(X_{\De_n}/\sqrt{\De_n}))$ by
$\E(f_j(\si U))$ for all $j\in J''$ it there is $s\in I\cap[0,1)$
such that the sequence $T_n\De_n^{(1-s)}$ is bounded.
\end{theo}

\section{Some technical tools}\label{S-PREL}
\setcounter{equation}{0}
\renewcommand{\theequation}{\thesection.\arabic{equation}}

\subsection{Estimates for $X''$.}

Let us begin with some additional notation. Below, $C$ denotes a
constant which may change from line to line and depend on the
characteristics $(b,c,F)$, and sometimes on the index $r$.

With any process $Y$ we associate its
discretized version $Y^{(n)}$ by
\bee\label{30}
Y^{(n)}_t~:=~Y_{\De_n[t/\De_n]}~=~\sum_{i=1}^{[t/\De_n]}\dd Y.
\eee
As soon as $Y$ is c\`adl\`ag, by a well known result
$Y^{(n)}(\om)\tosc Y$ (convergence for the Skorokhod topology).

Recalling the notation $X'$ and $X''$ of (\ref{14}), we have:

\begin{lem}\label{L1} We have $X''_t/\sqrt{t}\top0$ as
$t\to0$, and also  $X'_t/t\top0$ if $1\in I$, and
\bee\label{31}
\E(\phi_2(X''_t/\sqrt{t}))~=~
\left\{\begin{array}{ll}
\mbox{\rm o}\left(t^{1-s/2}\right)\qquad&\mbox{if }~ s\in
I\cap(0,2]\\[2mm]
\mbox{\rm O}\left(t\right)\qquad&\mbox{if }~ 0\in I.
\end{array}\right.
\eee
\end{lem}

\nib Proof. \rm If $s\in I$ with $s\leq2$ the
law of $X''_1$ belongs to the class $\g_s$ defined in the paper
\cite{AJ}. Hence $X''_t/\sqrt{t}\top0$, resp. $X'_t/t\top0$ if $1\in
I$, follows from Lemma 3.1--(a) of that paper applied
with $Y=X''$ and $\be=2$, resp. $Y=X'$ and $\be=1$.

In the same lemma, the estimate (\ref{31}) is proved when the process
$X$ is symmetrical and $s>0$, but a close look at the proof shows that it also
holds in general for $s>0$. Finally when $0\in I$ then $X''_t=\Bb
t+X'_t$ and $X'$ is compound Poisson, so (\ref{31}) is obvious. \qed

\subsection{Convergence of triangular arrays.}\label{SS2-1}

First, we gather in the following lemma some classical results (see
e.g. \cite{JS}) which will be heavily used below.
Let $(\chi^n_i=(\chi^{nj}_i)_{1\leq j\leq d})$ be a rowwise
i.i.d. array of $\R^d$--valued random variables, and
let $Z^n_t=\sum_{i=1}^{[t/\De_n]}\chi^n_i$ be their partial sums.

Let now $Z$ be a $d$--dimensional L\'evy process with characteristics
$(b',c',F')$: we have $b'=(b'^j)\in\R^d$, and
$c'=(c'^{ij})$ is a nonnegative symmetric matrix, and $F'$ is a positive
measure on $\R^d$. We also consider an arbitrary $C^1_b$
truncation function $\psi$ on $\R$: it satisfies $\psi(x)=x$ on a neighborhood
of $0$ and has compact support.

\begin{lem}\label{L4} a) If for some $t>0$ the sequence of
variables $Z^n_t$ converges in law to a limit, then the processes $Z^n$
converge in law, for the Skorokhod topology, to a L\'evy process.

b) We have $Z^n\tol Z$ if and only if for some (and then all) $C^1_b$
truncation function $\psi$ we have
\bea\label{Conv1}
\left.\begin{array}{l}
\frac1{\De_n}\E(\psi(\chi^{nj}_1))~\to~
b'^j+\int(\psi(x^j)-x^j1_{\{|x^j|\leq1\}})F'(dx),~~~1\leq j\leq d\\[2.5mm]
\frac1{\De_n}\E(\psi(\chi^{ni}_1)\psi(\chi^{nj}_1))~\to~
c'^{ij}+\int \psi(x^i)\psi(x^j)F'(dx),~~~1\leq i,j\leq d\\[2.5mm]
\frac1{\De_n}\E(f(\chi^n_1))~\to~\int f(x)F'(dx),~~~f~~
\mbox{bounded $F'$-a.e. continuous
  null around }~0.\end{array}\!\!\right\}
\eea
\end{lem}

Next we turn to the so--called {\em stable convergence in law}, which
has been defined in Subsection \ref{SN}.
We use this kind of convergence in the following setting: we are
in the situation of Lemma \ref{L4}--(b) with $Z=Z'+W'$, where $Z'$ is a
L\'evy process on the space $\proba$, with characteristics
$(b',0,F')$, and $W'$ is a $d$--dimensional Wiener process with unit
variance--covariance matrix $c'$, defined on an
extension of the original space and independent of $X$: so
$Z$ is a L\'evy process with characteristics $(b',c',F')$.

\begin{lem}\label{LL1} Let $Z=Z'+W'$ be as above, and suppose that
$\chi^n_i$ is $\f_{i\De_n}$--measurable and satisfies the first two
properties in (\ref{Conv1}) with $k_n=1/\De_n$. If
\bee\label{60}
\frac1{\De_n}~\PP(\|\chi^{n}_1-Z'_{\De_n}\|>\ep)\to0\qquad\forall\ep>0,
\eee
\bee\label{61}
\frac1{\De_n}~\E(\psi(\chi^{nj}_1) X^c_{\De_n})\to0
\eee
for all $j=1,\cdots,d$, the processes $Z^n$ converge stably in law to
$Z$, as defined above.

If further $c'=0$ (equivalently $W'=0$) (\ref{60}) alone implies
$Z^n\top Z=Z'$.
\end{lem}

\nib Proof. \rm Set $U^n=Z^n-Z'^{(n)}$, and suppose for a moment that
the following two properties hold:
\bee\label{62}
\left.\begin{array}{l}
\mbox{ under (\ref{60}) and
  (\ref{61}), the processes $U^n$ converge stably in law to $W'$ }\\[2mm]
\mbox{ under (\ref{60}) and if $c'=0$, the processes $U^n$
  converge in law to $0$ }.
\end{array}\right\}
\eee
We have $Z'^n\tosc Z'$
pointwise. Then the second part of (\ref{62}) gives our
second claim. The first part of (\ref{62}) yields that the sequence
$U^n$ is C-tight (that is, tight with all limiting processes being
continuous), so
obviously for any continuous bounded function $g$ on the
Skorokhod space we have $g(Z'^n+U^n)-g(Z'+U^n)\top0$.
If further $V$ is bounded $\f$--measurable, then
$\E(g(Z^n)V)-\E(g(Z'+U^n)V)\to0$. Applying once more
(\ref{62}) and the fact that the map $u\mapsto g(Z'(\om)+u)$
(where $u$ is a c\`adl\`ag function) is continuous for the Skorokhod
topology at any function $u$ which is continuous in time, and
remembering the definition of stable convergence, we see that
in fact $\E(g(Z^n)V)\to\E(g(Z'+W')V)$, which is our first claim.

So it remains to prove (\ref{62}). For this, we observe that
$U^n_t=\sum_{i=1}^{[t/\De_n]}\ze^n_i$, where the $\ze^n_i$'s
for $i=1,2,\ldots$ are i.i.d. and distributed as
$\chi^n_1-Z'_{\De_n}$. First, $Z'^{(n)}\tosc Z'$ as said before, and $Z'$
has characteristics $(b',0,F')$, hence Lemma \ref{L4} implies that the
variables $Z'_{\De_n}$ satisfy (\ref{Conv1}) with
the same $b'$ and $F'$, and with $c'=0$. Second, (\ref{60}) obviously
implies that
the last part of (\ref{Conv1}) is satisfied by the variables $\ze^n_1$
with $F'=0$. Then consider the three properties:
\bee\label{63}
\frac1{\De_n}~\E(\psi(\ze^{nj}_1))~\to~0,
\eee
\bee\label{64}
\frac1{\De_n}~\E(\psi(\ze^{ni}_1)\psi(\ze_1^{nj}))~\to~c'^{ij},
\eee
\bee\label{65}
\frac1{\De_n}~\E(\psi(\ze^{nj}_1)\psi(X^c_{\De_n}))~\to~0.
\eee
By Lemma \ref{L4} we have the second part of (\ref{62}) if
(\ref{60}), (\ref{63}) and (\ref{64}) with $c'=0$ holds, and by the
criterion given in Theorem IX.7.19 of \cite{JS} if we additionally
have (\ref{65}) then the first part of (\ref{62}) holds.

The properties of $\psi$ yield a
continuous bounded function $g$ vanishing around $0$, such that
for all $\ep>0$ small enough we have
$|\psi(x-x')-\psi(x)+\psi(x')|\leq C1_{\{|x-x'|>\ep\}}
+\ep g(x)$. Then (\ref{60}) and the third part of (\ref{Conv1})
for $Z'_{\De_n}$ yield
$$\limsup_n\frac1{\De_n}~\left(\E(\psi(\ze^{nj}_1))-\E(\psi(\chi^{nj}_1))
+\E(\psi(Z'^j_{\De_n}))\right)\leq ~\ep\int F'(dx)g(x^j).$$
$\ep$ is arbitrarily small, hence the above $\limsup$ actually
vanishes and, since
$\frac1{\De_n}~\E(\psi(Z'^j_{\De_n}))$ and $\frac1{\De_n}~
\E(\psi(\chi^{nj}_1))$ have the same limit, we deduce (\ref{63}).

In a similar way, if $\psi_\ga(x)=\ga\psi(x/\ga)$, we have for all
$\ga>0$ a continuous bounded function $g_\ga$ vanishing
around $0$, such that for all $\ep>0$ small enough we have
\bean
&&|\psi(x^i-x'^i)\psi(x^j-x'^j)
-\psi(x'^i)\psi(x^j)+\psi(x'^i)\psi(x'^j)|~\leq~
C1_{\{\|x-x'\|>\ep\}}
+\ep g_\ga(\|x\|)\\
&&\hskip3cm+|\psi_\ga(x'^i)|(|\psi_\ga(x^j)|+|\psi_\ga(x'^j)|)
+|\psi_\ga(x'^j)|(|\psi_\ga(x^i)|+|\psi_\ga(x'^i)|).
\eean
Recall that $\frac1{\De_n}~\E(\psi_\ga(Z'^j_{\De_n})^2)\to
\int \psi_\ga(x^j)^2F'(dx)$, which in turn goes to $0$ as
$\ga\to0$, whereas $\frac1{\De_n}~\E(\psi_\ga(\chi^{nj}_1)^2)\to
c'^{jj}+\int \psi_\ga(x^j)^2F'(dx)$, which stays bounded as
$\ga\to0$. Then using Cauchy--Schwarz and the same argument as for
proving (\ref{63}), we get by taking the $\limsup$ in $n$, then
letting $\ep\to0$, then letting $\ga\to0$:
$$\frac1{\De_n}~\left(\E(\psi(\ze^{nj}_i))\psi(\ze^{nj}_1))
-\E(\psi(\chi^{ni}_1)\psi(\chi^{nj}_1))
+\E(\psi(Z'^i_{\De_n})\psi(Z'^j_{\De_n}))\right)\to 0.$$
Then (\ref{64}) follows.

Finally we turn to (\ref{65}). We $|\psi(x-x')\psi(y)-\psi(x)y+\psi(x')y|
\leq C1_{\{|x-x'|>\ep\}}+\ep g(x)+C|y|1_{\{|y|>\ga\}}$, for some
$\ga>0$. So exactly as before, with in addition the fact that
$\E(|X^c_{\De_n}|1_{\{|X^c_{\De_n}|>\ga\}}\leq C\De_n^2$ for any
fixed $\ga$, we get
$$\frac1{\De_n}~\E(\psi(\ze^{nj}_1)\psi(X^c_{\De_n}))-
\frac1{\De_n}~\E(\psi(\chi^{nj}_1)X^c_{\De_n})
+\frac1{\De_n}~\E(\psi(Z'^j_{\De_n})X^c_{\De_n})~\to~0.$$
Since the last term in the left side above vanishes (because $X^c$ is
independent of $Z'$), we deduce (\ref{65}) from (\ref{61}), and the
proof is finished. \qed

\section{Theorem \ref{T3} and some consequences}\label{S-LLNWN}
\setcounter{equation}{0}
\renewcommand{\theequation}{\thesection.\arabic{equation}}

\subsection{Proof of Theorem \ref{T3}--(i,ii).}

Before proceeding to the proof, let us introduce a family
$(\psi_\eta)_{\eta>0}$ of $C^2$ functions satisfying
\bee\label{44}
1_{[-\eta,\eta]}(|x|\leq\psi_\eta(x)
\leq1_{[-2\eta,2\eta]}(|x|),\qquad |\psi'_\eta|\leq C/\eta,
\quad|\psi''_\eta|\leq C/\eta^2
\eee
(take for example $\psi_\eta(x)=\psi(x/\eta)$, where $\psi$ is $C^2$
with $1_{[-1,1]}(|x|)=\leq\psi(x)\leq1_{[-2,2]}(|x|)$ and bounded derivatives).
Next, the proof goes through a number of steps.
\vst

\ni\it Step 1. \rm (ii) is a particular case of (i), and if $f$
satisfies either one of the conditions in (i) the process $f\star\mu$
is well defined and c\`adl\`ag, hence $V(f)$ as well. In view of
the convergence $V(f)^{(n)}\tosc V(f)$, which holds $\om$--wise, it is
enough to prove the following:
\bee\label{50}
\BZ^n(f):=V^n(f)-V(f)^{(n)}~\toucp~0.
\eee
\vst

\ni\it Step 2: \rm Here we prove (\ref{50}) when
$f\in C^{0,F}$ vanishes on a neighborhood of $0$, say
$[-2\ep,2\ep]$. We introduce a series of notation to
be used later again. Let $S_1,S_2,\cdots$ be the
successive jump times of $X$ with absolute size bigger than $\ep$,
and $R_p=\De X_{S_p}$ (these depend on $\ep$, of course); we put
\bee\label{92}
X(\ep)_t=X_t-(x1_{\{|x|>\ep\}})*\mu_t=X_t-
\sum_{p:~S_p\leq t}R_p,
\eee
and $R'^n_p=\dd X(\ep)$ on the set
$\{S_p\in((i-1)\De_n,i\De_n]$. Finally if $T>0$ let
$\Om_n(T,\ep)$ be the set of all $\om$ such that each interval
$[0,T]\cap((i-1)\De_n,i\De_n]$ contains at most one $S_p(\om)$, and
such that $|\dd X(\ep)(\om)|\leq2\ep$ for all $i\leq T/\De_n$. Then
$\Om_n(T,\ep)\to\Om$ as $n\to\infty$.

Recalling $f(x)=0$ when $|x|\leq2\ep$, we see that on the
set $\Om_n(T,\ep)$ and for all $t\leq T$,
\bee\label{91}
V^n(f)_t-V(f)^{(n)}_t=\sum_{p:~S_p\leq\De_n[t/\De_n]}
(f(R_p+R'^n_p)-f(R_p)),
\eee
Since $f$ is $F$--a.e. continuous there
is a null set $N$ such that, if $\om\notin N$, then $f$ is
continuous at each point $R_p(\om)$, whereas $R'^n_p(\om)\to0$, so
$\sup_{t\leq T}|V^n(f)_t-V(f)^{(n)}_t|\to0$ when
$\om\notin N$. Hence
(\ref{50}) is obvious (we even have almost sure convergence).
\vst

\ni\it Step 3: \rm For each $r$ let $f_r$ be a given nonnegative function
with compact support and $C^2$ outside any neighborhood of
$0$ and satisfying $f(x)=|x|^r$ on some neighborhood of $0$.
Here we prove that it is enough to have the following property: (\ref{50})
holds for $f_r$, when $r=2$, or $r\in I\cap[1,2)$ and $c=0$, or
$r\in I\cap(0,1)$ and $c=\Bb=0$: we assume this below.

Let $f\in\ea''_r\cap C^{0,F}$ with $r$ as above. For
any $\eta>0$ the function $f(1-\psi_\eta)$ satisfies (\ref{50})
with $V(f(1-\psi_\eta))=(f(1-\psi_\eta))*\mu$ by
Step 2, and we have $|f\psi_\eta|\leq\ep_\eta f_r$ with
$\ep_\eta$ going to $0$ as $\eta\to0$. Then  $|\BZ^n(f\psi_\eta)|
\leq\ep_\eta(|\BZ^n(f_r)|+|V(f_r)^{(n)}|)$. Since $f_r$
satisfies (\ref{50}) by hypothesis, and $V(f_r)^{(n)}\to
V(f_r)$, by letting $\eta\to0$ we readily deduce (\ref{50})
with $V(f)=f*\mu$ for $f$.
This proves [a-1] and [a-3] and, since $f-f_r$
belongs to $\ea_r''$ as soon as $f\in\ea_r$, this also proves (b)
and (c).

Finally let $f\in\ea'_r\cap C^{0,F}$ in cases [a-2] or [a-4].
On the one hand $f(1-\psi_\eta)$ and $f_r(1-\psi_\eta)$ satisfy (\ref{50})
by Step 2, and $f_r$ does by hypothesis, hence $f_r\psi_\eta$ does as well.
On the other hand $|f\psi_\eta|\leq Cf_r\psi_\eta$ for some constant
$C$, hence $|\BZ^n(f\psi_\eta)|\leq |\BZ^n(f_r\psi_\eta)|+
V(f_r\psi_\eta)^{(n)}$. Since $V(f_r\psi_\eta)^{(n)}$
goes to $0$ locally uniformly in time
as $\eta\to0$, we get (\ref{50}) for $f$.
\vst

\ni\it Step 4: \rm We are left to prove (\ref{50}) when $f=f_r$ as given
in Step 3, with $r$ as above. Set
$$g(x,y)=f(x+y)-f(x)-f(y)-xf'(y)1_{\{|x|\leq1\}},\quad
k(x,y)=f(x+y)-f(x)-f(y)$$
with the convention $f'(0)=0$ if $r\leq1$ (otherwise, $f'(0)$ is the
derivative of $f$ at $0$, of course). If we apply the usual It\^o's
Formula if $r=2$ (so $f$ is $C^2_b$ over $\R$) and its
extension given in Theorem 3.1 of \cite{JJM} if $r<2$, we see that
$f(X)-V(f)=A+M$, where $A_t=\int_0^t\al(X_s)ds$ and $M$ is a locally
square--integrable martingale, all these given by
\bee\label{51}
\al(y)=\left\{\begin{array}{ll}
bf'(y)+\frac c2(f''(y)-f''(0))
+\int F(dx)g(x,y)&\mbox{if }~ r=2\\[2mm]
bf'(y)+\int F(dx)g(x,y)
\quad&\mbox{if }~1<r<2,~c=0\\[2mm]
\Bb f'(y)-|\Bb|+\int F(dx)k(x,y)
\quad&\mbox{if }~r=1,~c=0\\[2mm]
\int F(dx)k(x,y)&\mbox{if }~0<r<1,~c=\Bb=0,
\end{array}\right.
\eee
$$M_t=\left\{\begin{array}{ll}
\int_0^tf'(X_s)dX^c_s+k(x,X_-)\is_t\quad&\mbox{if }~r=2\\[2mm]
k(x,X_-)\is_t &\mbox{if }~r<2
\end{array}\right.$$
Moreover the predictable bracket of $M$, say $A'$, has the form
$A'_t=\int_0^t\al'(X_s)ds$, where
\bee\label{52}
\al'(y)=\left\{\begin{array}{ll}
cf'(y)^2+\int F(dx)k(x,y)^2\quad&\mbox{if }~r=2\\[2mm]
\int F(dx)k(x,y)^2
\quad&\mbox{if }~r<2.
\end{array}\right.
\eee
\vst

\ni\it Step 5: \rm Observe that $Z^n(f)_t$ is the sum of $[t/\De_n]$
i.i.d. summands, distributed as $\ze_n=f(X_{\De_n})-V(f)_{\De_n}=
A_{\De_n}+M_{\De_n}$. Since $\E(M_{\De_n})=0$, it is clear that
(\ref{50}) follows from the two properties
$$ \frac1{\De_n}~\E(|A_{\De_n}|)~\to~0,\qquad
\frac1{\De_n}~\E(M_{\De_n}^2)~\to~0.$$
But $\E(|A_{\De_n}|)\leq\int_0^{\De_n}\E(|\al(X_t)|dt$ and
$\E(M_{\De_n}^2)=\int_0^{\De_n}\E(\al'(X_t)|dt$. Then it
is enough to prove that, as $t\to0$
\bee\label{53}
|\al(y)|\leq C,\quad |\al'(y)|\leq C,\quad
\al(X_t)\top0,\quad \al'(X_t)\top0.
\eee
\vst

\ni \it Step 6: \rm Finally we prove (\ref{53}). Suppose first $r=2$,
so $f$ is $C^2_b$ and $f(0)=f'(0)=0$. Hence $|k(x,y)|\leq C\phi_1(x)$ and
$|g(x,y)|\leq C\phi_2(x)$ and
$\lim_{y\to0}k(x,y)=\lim_{y\to0}g(x,y)=0$. Then the first two
properties of (\ref{53}) are obvious, and by Lebesgue's theorem
$\al(y)\to0$ and $\al'(0)\to0$ as $y\to0$, hence we get the last two
properties as well.

Suppose next that $1<r<2$, so the function $f$ is $C^1_b$ and $f'$ is
H\"older with index $r-1$,
and $f(0)=f'(0)=0$, hence $|k(x,y)|\leq C\phi_1(x)$ and
$|g(x,y)|\leq C\phi_r(x)$, and again
$\lim_{y\to0}k(x,y)=\lim_{y\to0}g(x,y)=0$. Then we conclude
(\ref{53}) as above, upon using the assumption $r\in I$.

Finally if $r\leq1$ we have $|k(x,y)|\leq
C\phi_r(x)$, and again $k(x,y)\to0$ as $y\to0$. Since $r\in I$
we conclude as above if $r<1$, or if $r=1$ and $\Bb=0$. If$r=1$ and
$\Bb\neq 0$, in order to get
$\al(X_t)\top 0$ we need an additional argument: we
use Lemma \ref{L1}, which yields $X'_t/t\top 0$, so if $\Bb>0$, say,
then $\PP(X_t>0)\to1$ because $X_t=X'_t+\Bb t$, and we have $f(y)=1$ if
$y>0$, so the result follows. The same argument applies when $\Bb<0$,
using $f'(y)=-1$ if $y<0$. \qed
\vsc

\subsection{The behavior of $H_t(f)$.}

We draw some consequences of the convergence
$V^n(f)\top V(f)$ in Theorem \ref{T3}--(i), which a fortiori implies
$V^n(f)\tol V(f)$. Note that $V^n(f)_t$ is the sum of $[t/\De_n]$ i.i.d.
variables
distributed as $f(X_{\De_n})$. Hence the first part of (\ref{Conv1})
applied with a truncation $\psi$ having $\psi(x)=x$ as soon as
$|x|\leq\sup|f|$ when $f$ is bounded, plus the fact that the above
convergence takes place
for any sequence $\De_n\to0$, yield the following, for any
$f\in C^{0,F}$:
\bee\label{Conv2}
\frac1{t}~H_{t}(f)\to\left\{\begin{array}{ll}
c+F(f)&\mbox{if }~ f\in\ea_2^b\\[2.5mm]
|\Bb|+F(f)\quad&\mbox{if }
~ f\in\ea_1'^b\\[2.5mm]
F(f)&\mbox{if }\left\{\begin{array}{l}
f\in\ea''_2\\
f\in\ea'_r,~~1<r<2,~~r\in I,~~c=0\\
f\in\ea''_1,~~1\in I,~~c=0\\
f\in\ea'_r,~~r<1,~~r\in I,~~c=\Bb=0\end{array}\right.
\end{array}\right.
\eee

We deduce from these facts a number of useful properties.

\begin{lem}\label{L5} Assume $c=0$ and let $r\in(0,2)$. Then
\bee\label{Conv3}
H_{t}(\phi_r)=\left\{\begin{array}{ll}
\mbox{\rm o}(t^{r/s})\quad &\mbox{if $r<s$ with $s\in I$ and further
$\Bb=0$ when $s\leq1$}\\[2mm]
\mbox{\rm O}(t^{r}) &\mbox{if $r\leq1$ and $1\in I$ and
$\Bb\neq0$ when $s\leq1$},
\end{array}\right.
\eee
and for any $a>0$ and $s\in I\cap[0,2]$,
\bee\label{Conv4}
\E(\phi_r(X_t)1_{\{|X_t|\geq a\sqrt{t}\}})=\left\{\begin{array}{ll}
\mbox{\rm o}(t^{1-s/2+r/2})\quad &\mbox{if }~r<s~~\mbox{and either }
s\geq1,
\\&\mbox{or $s<1$ and $\Bb=0$}\\[2mm]
\mbox{\rm o}(t^{1-s/2+rs/2}) &\mbox{if
  }~r<1,~0<s\leq1,~\Bb\neq0\\[2mm]
\mbox{\rm O}(t) &\mbox{in all other cases.}
\end{array}\right.
\eee

\end{lem}

\nib Proof. \rm When $r\leq1$, H\"older
inequality yields $t^{-r}H_{t}(\phi_r)\leq \left(t^{-1}H_{t}(\phi_1)
\right)^r$, hence the second part of
(\ref{Conv2}) yields the second part of (\ref{Conv3}). Moreover if
$\phi_{r,\ep}(x)=|x|^r\bigwedge\ep$, for $\ep\in(0,1)$ and $r<s$
we have by H\"older inequality again
$$t^{-r/s}H_{t}(\phi_r)\leq
\left(t^{-1}H_{t}(\phi_{r,\ep}^{s/r})
\right)^{r/s}+t^{-r/s}H_{t}(\phi_r-\phi_{r,\ep}).$$
Since $\phi_r-\phi_{r,\ep}\in\ea_2^b\cap C^0$ the last term above goes to
$0$ for any $\ep$ because $r<s$, by the third part of (\ref{Conv2}).
Since $\phi_{r,\ep}^{s/r}\in\ea_s^b\cap C^0$, the first
term on the right goes to $F(\phi_{r,\ep}^{s/r})$, which in turn
goes to $0$ as $\ep\to0$ because $s\in I$. Then we obtain the first
part of (\ref{Conv3}).

Now if
$r<s$ and either $s\geq1$ or $s<1$ an $\Bb=0$, H\"older inequality yields
\bean
\E(\phi_r(X_t)1_{\{X_t\geq a\sqrt{t}~\!\!\}})
&\leq & \Big(\E(\phi_s(X_t))\Big)^{r/s}
\left(\PP(|X_t|\geq a\sqrt{t})\right)^{1-r/s}\\
&\leq&
C_aH_t(\phi_s)^{r/s}\left(\E(\phi_2(X_t/\sqrt{t}))\right)^{1-r/s}.
\eean
Then we apply (\ref{31}) (we have $X=X''$ here), and also
(\ref{Conv2}) which gives $H_t(\phi_s)\leq Ct$, and we get the
firts part of (\ref{Conv4}). When $s\leq 1$ and $\Bb\neq0$ another
application of H\"older inequality yields for $r\leq 1$:
$$\E(\phi_r(X_t)1_{\{X_t\geq a\sqrt{t}~\!\!\}})
\leq
C_aH_t(\phi_1)^{r}\left(\E(\phi_2(X_t/\sqrt{t}))\right)^{1-r}.$$
and we conclude the second part of (\ref{Conv4}), and also the third
part when further $s=0$. Finally if $r\geq s$ and further $\Bb=0$ if
$s\leq1$ and $r<1$, the left side of (\ref{Conv4}) is smaller than
$H_t(\phi_r)$, hence (\ref{Conv2}) gives the last part of
(\ref{Conv4}) in these cases. \qed
\vsc

\begin{lem}\label{L2} a) If $f\in\ea^b_r$ for some $r\in(0,2)$ and
  $c>0$, we have
\bee\label{33}
t^{-r/2}H_{t}(f)~\to~c^{r/2}\mu_r,\qquad
t^{-1/2-r/2}\E(X^c_tf(X_t))~\to~0.
\eee
and also, when $r<1$ and $s\in I\cap[0,2]$,
\bee\label{34}
t^{-r/2}H_{t}(f)-c^{r/2}\mu_r~=~
\left\{\begin{array}{ll}
\mbox{\rm O}\left(t^{1-r/2}\right)\quad&\mbox{if }~
r\geq s,~\Bb=0, \mbox{or if }~s=0\\[2mm]
\mbox{\rm o}\left(t^{1-s/2}\right)\quad&\mbox{if }~r<s
~~\mbox{and either $s\geq1$}\\ &\mbox{or $s<1$ and}~\Bb=0\\[2mm]
\mbox{\rm o}\left(t^{1-s/2-r/2+rs/2}\right)~~~&\mbox{if }~\Bb\neq0,~0<s\leq1.
\end{array}\right.
\eee
Moreover when there is no jump we have (\ref{33}) and also
$t^{-r/2}H_{t}(f)-c^{r/2}\mu_r=\mbox{\rm O}(t)$ for all $f\in\ea_r$
and all $r\in(0,\infty)$.

b) If $f$ is a bounded function we have
\bee\label{35}
\E(f(X_t/\sqrt{t}))~\to~\E(f(\si U)),\qquad
t^{-1/2}\E(X^c_tf(X_t/\sqrt{t}))~\to~0.
\eee
If further $f$ is an even function, then
\bee\label{36}
\E(f(X_t/\sqrt{t}))-\E(f(\si U))~=~
\left\{\begin{array}{ll}
\mbox{\rm o}\left(t^{1-s/2}\right)\quad&\mbox{if }s\in I\cap(0,2)\\[2mm]
\mbox{\rm O}(t)&\mbox{if }~0\in I.\end{array}\right.
\eee

c) With the assumptions of (a) and if $f'$ is a bounded function we have
\bee\label{36'}
t^{r/2}\E(f(X_t)f'(X_t/\sqrt{t}))~\to~c^{r/2}\mu_r\E(f'(\si U)).
\eee
\end{lem}

It would of course be possible to have estimates like (\ref{34}) when
$r\geq 1$ or $s>1$, but they are useless for us in this paper.
\vsq

\nib Proof. \rm  First we prove a simple auxiliary result. We denote
by $h$ the density of $\n(0,1)$, and $h_p(x)=x^ph(x)$ for $p\in\N$.
Let $r\geq0$ and $\al>1$ and $g$ be a nonnegative function such that
$g(x)\leq|x|^r$ if $r>0$ and $g(x)\leq1$ if $r=0$. For $p=0$ or $p=1$
we set
$$k_{p,\al}(y)=\E\left(U^p\left(g(y+U)\bigwedge\al^r\right)\right)
-\E(U^pg(U)).$$
By a change of variable we have
$$k_{p,\al}(y)=\int(h_p(x-y)-h_p(x))g(x)dx
-\int h_p(x-y)\left(g(x)-g(x)\bigwedge\al^r\right)dx.$$
The function $h_p$ is $C^\infty$, decreases exponentially
fast as $|x|\to\infty$, and its two first derivatives are
$h'_p(x)=x^{p-1}(p-x^2)h(x)$ and
$h''_p(x)=x^p(x^2-2p-1)h(x)$ (recall $p=0$ or $p=1$). We also have
$g(x)-g(x)\bigwedge\al^r\leq |x|^r1_{\{|x|>\al\}}$. Hence the
function $k_{p,\al}$ satisfies
$$|y|\leq1~~\Rightarrow~~
\left|k_{p,\al}(y)-y\int x^{p-1}(x^2-p)h(x)g(x)dx\right|\leq
C(y^2+\al^{-2}).$$
On the other hand $g(U+y)\leq C(|U|^r+\al^r\phi_r(y/\al))$, hence
$|k_{p,\al}(y)|\leq C(1+\al^r\phi_r(y/\al))$. Putting all these facts
together yiels:
\bee\label{37}
|k_{p,\al}(y)|\leq\left\{\begin{array}{l}
C\left(\al^{-2}+\phi_1(y)+\al^r\phi_r(y/\al)1_{\{|y|>1\}}\right)\\[2mm]
C\left(\al^{-2}+\phi_2(y)+\al^r\phi_r(y/\al)1_{\{|y|>1\}}\right)
\quad\mbox{if~ $x\mapsto x^pg(x)$ ~is even.}
\end{array}\right.
\eee

a) Let $f\in\ea_r^b$ with $r\in(0,2)$. Then $|f-\phi_r|\leq g$ for
some $g\in\ea_2''^b\cap C^0$, hence (\ref{Conv2}) yields
$\left|H_t(f)-H_t(\phi_r)\right|\leq Ct$. Also we
obviously have $|\E(X^c_t(f-\phi_r)(X_t))|\leq
\si\sqrt{t}~H_t((f-\phi_r)^2)^{1/2}\leq Ct^{3/2}$. Therefore it is
enough to prove (\ref{33}) and (\ref{34}) for $f=\phi_r$.

Set $\be_0(t)=t^{-r/2}H_{t}(\phi_r)-c^{r/2}\mu_r$ and
$\be_1(t)=t^{-r}\E(X^c_t\phi_r(X_t))$. Recall that the pair $(X_t,X^c_t)$ has
the same law as $(\si U\sqrt{t}+X''_t,\si U\sqrt{t})$, with $U$
independent of $X''_t$. Then we get for $p=0$ or $p=1$:
\bean
\be_p(t)&=&\si^{p+r}\left(\E\left(U^p~\left(|U+X''_t/\si\sqrt{t}|^r
\bigwedge(ct)^{-r/2}\right)\right)-\E(U^p|U|^r)\right)\\
&=&\si^{p+r} \E(k_{p,1/\si\sqrt{t}}(X''_t/\si\sqrt{t})),
\eean
where $k_{p,\al}$ is associated with function $g(x)=|x|^r$.
The first estimate in (\ref{37}) yields
$$|\be_p(t)|\leq C\left(t+\E(\phi_2(X''_t/\si\sqrt{t}~\!\!))^{1/2}
+t^{-r/2}\E(\phi_r(X''_t))\right).$$
Then $\be_p(t)\to0$, that is (\ref{33}), follows from (\ref{31}) and
from (\ref{Conv3}) applied to $X''$ and with $s=2$. Next, the second
estimate in (\ref{37}) yields, since $g$ is even:
\bee\label{37'}
|\be_0(t)|\leq C\left(t+\E(\phi_2(X''_t/\si\sqrt{t}~\!\!))
+t^{-r/2}\E(\phi_r(X''_t)1_{\{|X``_t|>\si\sqrt{t}~\!\!\}})\right).
\eee
Therefore (\ref{34}) readily follows from (\ref{31}) and (\ref{Conv4}).

Finally when there is no jump we have $X''_t=\Bb t$, so (\ref{37'}) gives
$|\be_0(t)|\leq Ct$ as soon as $|bt|\leq\si\sqrt{t}$, hence the final
claim of (a).

b) Now let $f$ be bounded. Exactly as above, we get for
$p=0$ or $p=1$:
$$\E\left((X^c_t)^pf(X_t/\sqrt{t})\right)
-t^{p/2}\E((\si U)^pf(\si U))=t^{p/2}c^{p/2}\E(k_{p,\infty}(X''/\sqrt{t})),
$$
where $k_{p,\infty}$ is associated with the bounded function
$g(x)=f(\si x)$ and $\al=\infty$. Here $k_{p,\infty}$ is bounded, so
instead of (\ref{37}) we have
$|k_{p,\infty}(y)|\leq C\phi_1(y)$ in general, and also
$|k_{0,\infty}(y)|\leq C\phi_2(y)$ if further $f$ is
even. Then we conclude (\ref{35}) and (\ref{36}) like in (a), by using
(\ref{31}).

c) This is proved exactly as (\ref{33}) of (a), except that we take
$g(x)=|x|^rf'(\si x)$. \qed

\begin{lem}\label{L3} Assume $c=0$ and $1\in I$, and let $f\in\ea^b_r$
  for some $r<1$. Then
\bee\label{38}
t^{-r}H_{t}(f)~\to~|\Bb|^r
\eee
as $t\to0$.
Moreover we have
\bee\label{39}
r\in I,~~ f\in C^{0,F}\quad\Rightarrow\quad
\frac1t\left(H_t(f)-t^r|\Bb|^r\right)~\to~F(f).
\eee
\bee\label{40}
2r\in I\quad\Rightarrow\quad
\frac1{\sqrt{t}}\left(H_t(f)-t^r|\Bb|^r\right)~\to~0.
\eee
\end{lem}

\nib Proof. \rm  Exactly as in the proof of (a) in the previous lemma,
it is enough to prove the results for $f=\phi_r$. If
$H'_t(f)=\E(f(X'_t))$, we can apply (\ref{Conv3}) to
the process $X'$, whose drift $\Bb'$ vanishes and whose associated set
$I'$ is equal to $I$. This gives
\bee\label{41}
\left.\begin{array}{ll}
H_t(\phi_r)=\mbox{o}(t^r)\qquad&\mbox{(apply (\ref{Conv3}) with $s=1$)}\\[2mm]
2r\in I~~\Rightarrow~~H_t(\phi_r)=\mbox{o}(t^{1/2})\qquad&
\mbox{(apply (\ref{Conv3}) with   $s=2r$)}
\end{array}\right\}
\eee
Since $X_t=\Bb t+X'_t$ and
$\Big||x+y|^r-|x|^r\Big|\leq|y|^r$, it is easy to check that
$A_t=\phi_r(X_t)-|\Bb t|^r|-\phi_r(X'_t)$
satisfies $|A_t|\leq2|\Bb t|^r$ and also $|A_t|\leq2\phi_r(X'_t)$.

First we deduce $|\phi_r(X_t)-|\Bb t|^r|\leq 3\phi_r(X'_t)$,
hence $\left|H_{t}(\phi_r)-t^r|\Bb|^r\right|\leq 3H'_t(\phi_r)$, and
(\ref{38}) and (\ref{40}) follow from (\ref{41}).

Next, with the notation (\ref{44}) we have $|A_t|\leq
2(\phi_r\psi_\eta)(X'_t)+2|\Bb t|^r
1_{\{|X'_t|>\eta\}}$, and thus also
$|A_t|\leq 2(\phi_r\psi_\eta)(X'_t)+C_\eta t^r\phi_2(X'_t)$. Therefore
$$\left|H_t(\phi_r)-t^r|\Bb|^r-H'_t(\phi_r)\right|\leq
2H'_t(\phi_r\psi_\eta)+C_\eta t^rH'_t(\phi_2).$$
If $r\in I$, then (\ref{Conv2}) applied to $X'$ yields that
$\frac1t~\! H'_t(\phi_r)\to F(\phi_r)$, and also
$\frac1t~\! H'_t(\phi_r\psi_\eta)\to F(\phi_r\psi_\eta)$ and
$\frac1t~\! H'_t(\phi_2)\to F(\phi_2)$. We deduce that
$$\limsup\left|\frac1t~(H_t(\phi_r)-t^r|\Bb|^r)F(\phi_r)\right|
\leq 2F(\phi_r\psi_\eta)$$
and, since $F(\phi_r\psi_\eta)\to0$ as $\eta\to0$, we finally deduce
(\ref{39}). \qed

\begin{lem}\label{L6} Either we are in one of the cases described
in Theorem \ref{T3}--(ii) and every $f\in\ea_r^b$ satisfies
$\sup_t\frac1t~\! |H_{t}(f)|<\infty$, or $\frac1t~ H_{t}(\psi_r)\to\infty$ and
$H_t(f)/H_t(\psi_r)\to1$ for all $f\in\ea_r^b$.
\end{lem}

\nib Proof. \rm Everything is a simple consequence of (\ref{Conv2}),
(\ref{33}) and (\ref{38})). \qed
\vsc

\begin{lem}\label{L12}  a) If either $f\in\ea''^b_1$, or $f\in\ea'^b_r$
  for some $r<1$ with $2r\in I$ and $c=0$, the family $\frac1t~\!\Ga_t(f)$
  is bounded, and it converges to $F(f^2)$ if further $f\in C^{0,F}$.

b) If $r=1$ the family $\frac1t~\Ga_t(f)$ is bounded for all
$f\in\ea_r^b$, and it converges to $F(f^2)+c(1-\mu_1^2)$ if further
$f\in C^{0,F}$.

c) In all other cases (that is $2r\notin I$, or $r<1$ and $c>0$), we
have $\frac1t~\!\Ga_t(\phi_r)\to\infty$ and
$\Ga_t(f)/\Ga_t(\phi_r')\to1$ for all $f\in\ea_r^b$.
\end{lem}

\nib Proof. \rm (b) and the case $f\in\ea''^b_1$ of (a) and the case
$c>0$ of (c) are trivial consequences of Lemma \ref{L2} and of (\ref{Conv2}).

It remains to study the case $r<1$ and $c=0$. Let $f\in\ea_r'^b$.
Let also $p=1$ or $p=2$. We have $f^p-\phi_r^p\in\ea''^b_2$,
hence the family $\frac1t~\! H_t(f^p-\phi_r^p)$ is bounded, and converges to
$F(f^p-\phi_r^p)$ if further $f$ is $F$--a.e. continuous (apply
(\ref{Conv2})). Moreover (\ref{Conv3}) applied with $s=2$ yields
$H_t(f+\phi_r)\leq CH_t(\phi_r)=$ o$(t^{r/2})$. Since
$$\frac1t~(\Ga_t(f)-\Ga_t(\phi_r))=\frac1t~ H_t(f^2-\phi_r^2)
-\frac1t~H_t(f-\phi_r)H_t(f+\phi_r),$$
we deduce that the family $\frac1t~\:\!(\Ga_t(f)-\Ga_t(\phi_r))$ is
bounded, and converges to $F(f^2-\phi_r^2)$ if further $f$ is
$F$--a.e. continuous. We deduce the following properties:
\bee\label{43}
\left.\begin{array}{lll}
\frac1t~\Ga_t(\phi_r)\to F(\phi_r^2)<\infty
&\Rightarrow& \left\{\begin{array}{ll}
\frac1t~\Ga_t(f)\to F(f^2)&\forall f\in\ea_r'^b\cap C^{0,F}\\
\sup_t\left|\frac1t~\Ga_t(f)\right|<\infty~&\forall f\in\ea_r'^{b}
\end{array}\right.\\[4mm]
\frac1t~\Ga_t(\phi_r)\to\infty
~&\Rightarrow& \left\{\begin{array}{ll}
\frac1t~\Ga_t(f)\to\infty &\forall f\in\ea_r'^{b}\\
\Ga_t(f)/\Ga_t(\phi_r)\to1~~&\forall f\in\ea_r'^{b}.
\end{array}\right.
\end{array}\!\!\!\right\}
\eee

Suppose now that $2r\notin I$. Define $\psi_\eta$ by (\ref{44}), so that
$\frac1t~\!(\Ga_t(\phi_r)-\Ga_t(\phi_r\psi_\eta))\to
F(\phi_r^2-(\phi_r\psi_\eta)^2)$ by what precedes. Since
and $\Ga_t(\phi_r\psi_\eta)\geq0$, we get
$$\liminf_{t\to0}\frac1t~\Ga_t(\phi_r)~\geq~F(\phi_r^2(1-\psi_\eta^2)),$$
and since $F(\phi_r^2)=\infty$ by hypothesis and $1-\psi_\eta^2$
increases to $1$ as $\eta\to0$, we obtain
$\frac1t~\Ga_t(\phi_r)\to\infty$: this prove (c) in the case $c=0$.

Finally we suppose that $r<1$ and $2r\in I$ and $c=0$. In order to get
the remaining part of (a) it suffices to show
$\frac1t~\Ga_t(\phi_r)\to F(\phi_r^2)$. We
single out the following three cases:

1) Assume either $r>1/2$, or $r\leq1/2$ and $\Bb=0$. Then
(\ref{Conv2}) yields $\frac1t~ H_t(\phi_r^2)\to F(\phi_r^2)$, and (\ref{Conv3})
   applied with $s=2r$ yields $H_t(\phi_r)=$ o$(\sqrt{t})$, hence the
   result is obvious.

2) Assume $r=1/2$ and $\Bb\neq0$. (\ref{Conv3}) yields
$\frac1t~ H_t(\phi_r^2)\to F(\phi_r^2)+|\Bb|$, and (\ref{38}) gives
$\frac1{\sqrt{t}}~H_t(\phi_r)\to\sqrt{|\Bb|}$, so we obtain the result.

3) Assume $r<1/2$ and $\Bb\neq0$. Then (\ref{39}) and (\ref{40}) yield
$$\frac1t~H_t(\phi_r^2)=t^{2r-1}|\Bb|^{2r}+F(\phi_r^2)+\mbox{o}(1),\quad
\frac1{\sqrt{t}}~H_t(\phi_r)=t^{r-1/2}|\Bb|^{r}+\mbox{o}(1),$$
and the result follows. \qed

\subsection{Proof of Theorem \ref{T3}--(iii).}

We suppose that we are not in the cases of (ii) of this theorem. The
first claim of (iii) follows from Lemma \ref{L6}. For the second claim, it
is enough to prove it when $f=\phi_r$, because otherwise
$|f-\phi_r|\leq g$ for some $g\in\ea_2''\cap C^0$ and we know that
$\frac{\De_n}{H_{\De_n}(\phi_r)}~V^n(h)\toucp0$ by (i) of the theorem and
$H_{\De_n}(\phi_r)/\De_n\to\infty$.

So it is clearly enough to prove that $Y^n_t:=
\frac{\De_n}{H_{\De_n}(\phi_r)}~V^n(f)_t- \De_n[t/\De_n]\toucp0$.
Note that $Y^n_t$ is the sum of $[t/\De_n]$
i.i.d. variables distributed as $\ze_n=\De_n(\phi_r(\dd
X)/H_{\De_n}(\phi_r)-1)$. Then $\E(\ze_n)=0$ and
$$\E(\ze_n^2)=\frac{\De_n^2\Ga_{\De_n}(\phi_r)}{H_{\De_n}(\phi_r)^2}~
\leq \frac{\De_n^2H_{\De_n}(\phi_r^2)}{H_{\De_n}(\phi_r)^2}~
\leq C\frac{\De_n^2}{H_{\De_n}(\phi_r)}$$
because $\phi_r\leq1$. Hence $E(\ze_n^2)/\De_n\to0$, and
$Y^n\toucp0$ readily follows. \qed

\section{Proofs of the other theorems about $V^n(f)$ and $V'^n(f)$}\label{POT}
\setcounter{equation}{0}
\renewcommand{\theequation}{\thesection.\arabic{equation}}

\subsection{Proof of Theorem \ref{T4}.}

(i) It is enough to prove that the processes
$Y^n_t=\De_n^{1-r/2}V^n(f)_t-\De_nc^{r/2}\mu_r[t/\De_n]$
converge u.c.p. to $0$. Now $Y^n_t$ is the sum of $[t/\De_n]$
i.i.d. variables, all distributed as $\ze_n=\De_n(\De_n^{-r/2}f(\dd
X)-c^{r/2}\mu_r)$. If $f\in\ea_r^b$ and $r\in(0,2)$, or if $f\in\ea_r$
and $r\in[2,\infty)$ when there is no jump, we have
$$\E(\ze_n)=\De_n(\De_n^{-r/2}H_{\De_n}(f)-c^{r/2}\mu_r),\quad
\E(\ze_n^2)\leq \De_n^{2-r}H_{\De_n}(f^2).$$
Then (\ref{33}), and also (\ref{Conv2}) when $r\geq1$, yield
$E(\ze_n)/\De_n\to0$ and $E(\ze_n^2)/\De_n\to0$, and
$Y^n\toucp0$ readily follows.

It remains to study the case $f\in\ea_r\backslash\ea_r^b$ when $r<2$.
We do like for Theorem \ref{T3}--(iii): we have
$|f-g|\leq h$ for some $g\in\ea_r^b$ and $h\in\ea''_2\cap C^0$,
and Theorem \ref{T3}--(i) yields $\De_n^{1-r/2}V^n(h)\toucp0$
(recall $r<2$); since (i) holds for $g$ it also holds for $f$.
\vst

(ii) Let $f$ be bounded. We do as above with
$Y^n_t=\De_nV'^n(f)_t-\De_n\E(f(\si U))[t/\De_n]$
and $\ze_n=\De_n(f(\dd X/\sqrt{\De_n})-\E(f(\si U)))$. We have
$$\E(\ze_n)=\De_n\left(\E(f(\dd X/\sqrt{\De_n}))-\E(f(\si U))\right),\quad
\E(\ze_n^2)\leq \De_n^2\E(f(\dd X/\sqrt{\De_n})^2).$$
Then (\ref{35}) yields $|E(\ze_n)|/\De_n\to0$ and $E(\ze_n^2)/\De_n\to0$, and
$Y^n\toucp0$ follows. \qed

\subsection{Proof of Theorems \ref{T5} and \ref{T6}--(ii).}

\ni\it Step 1. \rm Let $f\in\ea_r\cap C^{0,F}$ with either $1\leq r<2$, or
$0<r<1$ and $2r\in I$ and $c=0$. Let $\vp$ be
continuous with $1_{[-a,a]}\leq\vp\leq1_{[-a',a']}$ for some
$0<a<a'<\infty$, or $\vp\equiv1$ when $f$ is bounded, so
the process $\Si(f,\vp)$ is well defined. For simplicity we always
write $v_n=H_{\De_n}(f\vp)$. We want to prove that $V^n(f)_t-v_nt$
converges stably in law to $Z=\Si(f,\vp)+W'$ with $W'$ a Wiener
process independent of $X$ and with unit variance $c(\mu_2-\mu_1^2)$ when
$r=1$ and $c>0$, and $W'=0$ otherwise (hence we get $V^n(f)_t-v_nt\top
\Si(f,\vp)$).

Suppose for a moment that the result holds whenever $f$ is
bounded. Then if we start with an unbounded $f$, we set
$f_\eta=f\psi_\eta$ (recall (\ref{44})).
Clearly $f_\eta\vp=f\vp$ as soon as $\eta>a'$, so for such $\eta$ we
have on the one hand that $V^n(f_\eta)_t-v_nt$ converges stably in law
to $\Si(f_\eta,\vp)+W'$, and on the other hand $V^n(f-f_\eta)\top
(f-f_\eta)\star\mu$ by Theorem \ref{T3}--(i). Since further
$\Si(f-f_\eta,\vp)$
and $(f-f_\eta)\star\mu$ both go to $0$ locally uniformly (pointwise in
$\om$), it is clear that $V^n(f)_t-v_nt$ converges stably in law to
$\Si(f,\vp)+W'$.
\vst

\ni\it Step 2. \rm In the sequel we suppose $f$ bounded.
Observe that $v_n=\E((f\vp)(X_{\De_n}))\to0$ because $f\vp$ is bounded
and $(f\vp)(x)\to0$ as $x\to0$ and $X_{\De_n}\to0$. Hence it is enough
to show the convergence of the processes $Z^n_t=V^n(f)_t-
v_n[t/\De_n]$ to $Z'=\Si(f,\vp)$.

Now, $Z^n_t$ is the sum of $[t/\De_n]$ i.i.d. variables
distributed as $\chi_n=f(X_{\De_n})-v_n$. Also, the
characteristics of $Z'$ are $(b',0,F')$ where
$b'=\int f(x)(1_{\{|f(x)|\leq1\}}-\vp(x))F(dx)$
and $F'$ is the image of $F$ by $f$. Then in view of Lemma \ref{LL1}
it is enough to prove the following properties, for some $C^1_b$
truncation function
$\psi$ such that $\psi(x)=x$ when $|x|\leq 2\sup_y|f(y)|$ and with
$c'=0$ when $r\neq1$ and $c'=c(1-\mu_1^2)$ when $r=1$:
\bee\label{70}
\frac1{\De_n}~\E(\psi(\chi_n))~\to~F(f(1-\vp)),
\eee
\bee\label{71}
\frac1{\De_n}~\E(\psi(\chi_n)^2)~\to~c'+F(f^2),
\eee
\bee\label{72}
\frac1{\De_n}~\E(|\chi_n-Z'_{\De_n}|>\ep)~\to~0\qquad\forall\ep>0,
\eee
plus the following, when $r=1$ and $c>0$:
\bee\label{73}
\frac1{\De_n}~\E(\psi(\chi_n)X^c_{\De_n})~\to~0.
\eee

\ni\it Step 3. \rm By our choice of $\psi$, we have $\psi(\chi_n)=
f(X_{\De_n})-v_n$. Then the left side of (\ref{70}) is
$\frac1{\De_n}H_{\De_n}(f(1-\vp))$. Since $f(1-\vp)\in\ea''^b_2\cap
C^{0,F}_r$,
(\ref{Conv2}) yields (\ref{70}). Also, the left side of (\ref{71})
is $\frac1{\De_n}\left(\Ga_{\De_n}(f)+H_{\De_n}(f(1-\vp))^2\right)$,
which by what precedes has the same limit as $\frac1{\De_n}\Ga_{\De_n}(f)$.
Then Lemma \ref{L12}--(a,b) gives (\ref{71}). Finally, the left side of
(\ref{73}) is $\frac1{\De_n}\E(X^c_{\De_n}f(X_{\De_n}))$, hence
(\ref{73}) when $r=1$ and $c>0$ follows from (\ref{33}).
\vst

\ni\it Step 4. \rm It remains to prove (\ref{72}). Let $\ep\in(0,1)$ be such
that $f(x)=|x|^r$ for $|x|\leq \ep$, and use the notation
$X(\ep)$, $S_p$, $R_p$ given around (\ref{92}), and
also $Z''_t=Z'_t-\sum_{s\leq t}\De Z'_s1_{\{|\De X_s|>\ep\}}$.
We divide the sample space into three sets: $B_n=\{S_1>\De_n\}$,
$B'_n=\{S_1\leq\De_n<S_2\}$, and the complement
$B''_n=\{S_2\leq\De_n\}$. Observe that
$$\chi_n-Z'_{\De_n}=\left\{\begin{array}{ll}
f(X(\ep)_{\De_n})-v_n-Z''_{\De_n}~~&\mbox{on }~B_n\\[2mm]
f(X(\ep)_{\De_n}+R_1)-v_n-
Z''_{\De_n}-f(R_1)~~&\mbox{on  }~B'_n\end{array}\right.$$
and $B'_n$ is independent of $(X(\ep),Z'',R_1)$ and
$\PP(B'_n)=$ O$(\De_n)$ and $\PP(B''_n)=$ O$(\De_n^2)$, so
\bean
&&\frac1{\De_n}~\PP(|\chi_n-Z'_{\De_n}|>3\ep^r)~\leq~
\frac1{\De_n}~\PP\left(\left|f(X'(\ep)_{\De_n})-v_n-
Z''_{\De_n}\right|>3\ep^r\right)\\
&&\qquad+C
\PP\left(\left|f(X(\ep)_{\De_n}+R_1)-v_n-
Z''_{\De_n}-f(R_1)\right|>3\ep^r\right)+C\De_n.
\eean
Since $X(\ep)$ (resp. $Z''$) is a
L\'evy process whose L\'evy measure charges only the interval
$(-\ep,\ep)$ (resp. $(0,\ep^r)$), (\ref{Conv2}) yields
$\frac1{\De_n}~\PP(|X(\ep)_{\De_n}|>\ep)\to0$ and
$\frac1{\De_n}~\PP(|Z''_{\De_n}|>\ep^r)\to0$, whereas $v_n\to0$: hence
the first term on the right side above goes to $0$. We also
have $X(\ep)_{\De_n}\to0$ and $Z''_{\De_n}\to0$, and for almost all $\om$
the function $f$ is continuous at $R_1(\om)$,
hence the second term on the
right side above goes to $0$, as well as the third one of course.
Then (\ref{72}) holds, and we are finished. \qed

\subsection{Proof of Theorem \ref{T6}--(i).}

Let $f\in\ea_r$ for some $r<1$, and assume $c=0$ and $2r\notin I$.
The Lemma \ref{L12}--(c) yields that
$u_n=\De_n/\sqrt{\Ga_{\De_n}(\phi_r)}$ goes to $0$.

There is $h\in\ea_2''\cap C^0$ such that $|f-\phi_r|\leq h$, and
$|V^n(f)-V^n(\phi_r)|\leq V^n(h)$ and by Theorem \ref{T3}--(i) we know that
$V^n(h)\top\Si(h)$. Since $u_n\to0$, it is then obvious
that we only need to prove the result for $f=\phi_r$. We also have
$|u_nH_{\De_n}(\phi_r)|\leq u_n\to0$, so
it is enough to prove the
convergence of the processes $Y^n_t=u_n\left(V^n(f)_t-
H_{\De_n}(f)[t/\De_n]\right)$.

Note that $Y^n_t$ is the sum of $[t/\De_n]$
i.i.d. variables distributed as
$\chi_n=u_n(f(X_{\De_n})-H_{\De_n}(f))$. Moreover $X^c=0$, and
$|\chi_n|\leq u_n\to0$, so in view of Lemma \ref{LL1} it is
enough to prove that $\frac1{\De_n}\E(\chi_n)\to0$ and
$\frac1{\De_n}\E(\chi_n^2)\to t$. But $\Ga_{\De_n}(f)$ is the variance
of $f(X_{\De_n})$, so $\E(\chi_n)=0$ and
$\E(\chi_n^2)=\De_n$, and we are done. \qed

\subsection{Proof of Theorem \ref{T7}.}

(i--1) We first prove the result when
$f\in\ea_r$ for some $r<1$. There is $g\in\ea''_2\cap C^0$ such
that $|f-\phi_r|\leq g$, and Theorem \ref{T3}--(i) implies that
$V^n(g)$ converges to a finite--valued limiting process, so because
$r<1$ it is clearly enough to prove the result for $f=\phi_r$. Since
$\De_n^{1/2-r/2}H_{\De_n}(\phi_r)\to0$ by (\ref{33}), we can replace $t$
by $\De_n[t/\De_n]$ in the centering term of our pre--limiting
processes. Therefore we are left to proving the stable convergence in
law of processes which at time $t$ are sums of $[t/\De_n]$
i.i.d. variables distributed as $\chi_n=u_n(\ze_n-\E(\ze_n))$, where
$u_n=\De_n^{1/2-r/2}$ and $\ze_n=\phi_r(X_{\De_n})$.

Note that $u_n\to0$ and $\phi_r$ is bounded. Then, choosing
an adequate truncation function $\psi$, and since $\chi_n$ has mean
$0$, we will deduce the result from Lemma \ref{LL1}, provided we prove
\bee\label{74}
\frac1{\De_n}~u_n^2\left(\E(\ze_n^2)-\E(\ze_n)^2\right)
~\to~c^r(\mu_{2r}-\mu_r^2),
\eee
\bee\label{75}
\frac1{\De_n}~u_n\E(X^c_{\De_n}\ze_n)~\to~0.
\eee
These two properties follow from (\ref{33}) (because $2r<2$ here).

Now we turn to the case without jump, with $f\in\ea_r$ and
$r\geq1$. In view of the last statement of Lemma \ref{L2}--(a), we
obviously can replace in the centering term
$t\De_n^{-r/2}H_{\De_n}(\phi_r)$ by $tc^{r/2}\mu_r$, and also by
$[t/\De_n]\De_n^{1-r/2}H_{\De_n}(f)$. With this last version for
centering, we are exactly in the same situation than above, with $u_n$
and $\chi_n$ similar, but $\ze_n=f(X_{\De_n})$. We have to prove
(\ref{74}) and (\ref{75}), which again hold because of
(\ref{33}). However, there is a difference here, namely $u_n$ no
longer tends to $0$. So additionally to the previous properties
we have to prove the following Lindeberg condition:
$$\frac1{\De_n}~u_n^2\E(\ze_n^21_{\{|u_n\ze_n|>\ep\}})~\to~0
\qquad\forall\ep>0.$$
Since $|\ze_n|\leq C(\De_n^{r}+|X^c_{\De_n}|^{r})$
has the same law as $C(\De_n^{r}+\De_n^{r/2}|U|^{r})$,
this Lindeberg condition is obvious and we have the desired
convergence.
\vst

(ii--1) Since $f$ is bounded, we clearly can replace $t$
by $\De_n[t/\De_n]$ in the centering term here again. So we are in the
same situation as in the first part of (i), with $u_n=\sqrt{\De_n}$
and $\ze_n=f(X_{\De_n}/\sqrt{\De_n})$, and the only things to prove
are (\ref{75}), and also (\ref{74}) with $\E(f(\si U)^2)-\E(f(\si
U))^2$ on the right side. These two facts follow from
(\ref{35}).
\vst

Finally, the claims in (i--2) and (ii--2) readily follow,
upon using (\ref{34}) and
(\ref{36}) (observe that the right side of (\ref{34}) is always
~o$(\sqrt{t}~\!\!)$ when $r<1$). \qed

\subsection{Proof of Theorem \ref{T8}.}

\ni\it Step 1. \rm Here we prove (a). We denote by $(\g_t)$ the
filtration on the extended space generated by $X$ and $W'$ and all processes
$U_p1_{\{T_p\leq t\}}$.
When $g$ vanishes around $0$, then $Z(g)_t$ is well
defined because the associated sum in (\ref{90}) is finite, and
if further $F(g^2)<\infty$ then $Z(g)$ is clearly a locally square--integrable
$(\g_t)$--martingale whose predictable bracket is $cF(g^2)t$.

Now let $g\in\ea'_1$. We write $g_n=(1-\psi_{1/n})g$ (recall
(\ref{44})). By what precedes $Z(g_n)$ is well defined, and also
$$\E\left((Z(g_{1/n}-g_{1/m})_t)^2\right)=cF((g_{1/n}-g_{1/m})^2)t
\leq Ct\int_{\{|x|\leq 2/n\}}x^2F(dx)$$
for $m\geq n$, because $|g(x)|\leq C|x|$ near $0$. Then the
sequence $Z(g_1-g_{1/n})_t$ converges in $\LL^2$ to a limit naturally
denoted as $Z(g_1)_t$, and the process $Z(g_1)$ is a locally
square--integrable $(\g_t)$--martingale. It remains to put
$Z(g)=Z(g_1)+Z(g_1)$.
\vst

\ni\it Step 2. \rm Now we prove (b) when $f$ is $C^1$ and vanishes
on the interval $[-2\ep,2\ep]$ for some $\ep>0$. We use the notation
$X(\ep)$, $S_p$, $R_p$,
$R'^n_p$ and $\Om_n(T,\ep)$ given around (\ref{92}). We have
$\{S_p<\infty\}=\cup_q\{S_p=T_q<\infty\}$, so if we set $U'_p=U_q$ on
the set $\{S_p=T_q<\infty\}$ the process
$$Z'_t=\sum_{p:~S_p\leq t}\si f'(R_p)U'_p$$
has clearly the same distribution, conditional on the $\si$--field
$\f$, as $Z(f')$ given by (\ref{90}). So the claim amounts to
the stable convergence in law, toward $Z'$, for the sequence of
processes $Z^n(f)=\frac1{\sqrt{\De_n}}~\!\left(V^n(f)-V(f)^{(n)}\right)$.

On the set $\Om_n(T,\ep)$ we have (\ref{91}) for $t\leq T$, hence also
\bee\label{96}
Z^n(f)_t~=
~\sum_{p:~S_p\leq\De_n[t/\De_n]}f'(R_p+\widetilde{R}'^n_p)
\frac{R'^n_p}{\sqrt{\De_n}}
\eee
for all $t\leq T$,
where $\widetilde{R}^n_p$ is between $R_p$ and $R_p+R'^n_p$.
Since $R^n_p\to 0$ and $\Om_n(T,\ep)\to\Om$, the result will easily
follow from the following property, for any integer $q\geq1$:
\bee\label{93}
\mbox{the random vectors $(R'^n_p/\sqrt{\De_n})_{1\leq
    p\leq q}$ converge stably in law to $(\si U'_p)_{1\leq p\leq q}$}.
\eee
We have $R'^n_p/\sqrt{\De_n}=A^n_p+A'^n_p$, where $A^n_p=\dd
X^c/\sqrt{\De_n}$ and $A'^n_p=\dd (X(\ep)-X^c)/\sqrt{\De_n}$ on the
set $\{S_p\in((i-1)\De_n,i\De_n]$. Observe that, since $X(\ep)$
is independent of $S_p$, the variable $A'^n_p$ has the same law as
$(X(\ep)-X^c)_{\De_n}/\sqrt{\De_n}$, which goes to $0$ in probability
by Lemma \ref{L1}. Hence it is enough to prove (\ref{93}) when
$R'^n_p$ is substituted with $A^n_p$.

On the other hand, the random vector $(A^n_p)_{1\leq p\leq q}$ is
distributed, conditionally on $\Om_n(T,\ep)$, as the vector
$(\si U'_p)_{1\leq p\leq q}$, and it is clearly asymptotically
independent of the process $X^c$, hence of the process $X$ as well.
Then we have (\ref{93}) and the result is proved.
\vst

\ni\it Step 3. \rm  Here we prove (b) for $f$ of class $C^1$, and
twice differentiable in a neighborhood of $0$ and $f''(x)=$ o$(|x|)$
as $x\to0$.
Set $f_\eta=f\psi_\eta$. Step 2 yields that $Z^n(f-f_\eta)$ converges
stably in law to $Z(f(1-\psi_\eta))$ for any $\eta>0$, whereas the
argument of Step 1 yields that $Z(f(1-\psi_\eta))\toucp Z(f)$ as
$\eta\to0$. Then the result follows from the property
\bee\label{94}
\lim_{\eta\to0}~\limsup_n~\PP\left(\sup_{t\leq T}
\left|Z(f_\eta)^n_t\right|>\ep\right)=0,\quad\forall\ep>0,~\forall T>0.
\eee

We know that $Z(f_\eta)^n_t$ is the sum of $[t/\De_n]$ i.i.d. variables
distributed as
$\ze(\eta)_n=(f_\eta(X_{\De_n})-V(f_\eta)_{\De_n})/\sqrt{\De_n}$.
Furthermore, the assumptions on $f$ imply that for $\eta$ small enough
then $f_\eta$ is $C^2$. Then if we
look at Step 4 of the proof of Theorem \ref{T3}--(i) we have
$f_\eta(X)-V(f_\eta)=A(\eta)+M(\eta)$ where $A(\eta)$ and $M(\eta)$ are
associated with $f_\eta$ as $A$ and $M$ are with $f$. Since
$\E(M(\eta)_t)=0$, (\ref{94}) will obviously follow from the two
properties
$$\lim_{\eta\to0}~\limsup_n~\De_n^{-3/2}\E(|A(\eta)_{\De_n}|)=0,\qquad
\lim_{\eta\to0}~\limsup_n~\De_n^{-2}\E(M(\eta)_{\De_n}^2)=0.$$
Now we associate the functions $\al_\eta$ and $\al'_\eta$ with $f_\eta$
by (\ref{51}) and (\ref{52}) (use the versions for $r=2$ in these
formula, and note that $f''_\eta(0)=0$; note also that $g=g_\eta$ and
$k=k_\eta$ depend on $\eta$ here). Recalling that $\E(|A(\eta)_t|)\leq
\int_0^tH_s(|\al_\eta|)ds$ and
$\E(M(\eta)_t^2)=\int_0^tH_s(\al'_\eta)ds$, it thus remains to prove
that
\bee\label{95}
\lim_{\eta\to0}~\limsup_{t\to0}~t^{-1/2}H_t(|\al_\eta|)=0,\qquad
\lim_{\eta\to0}~\limsup_{t\to0}~t^{-1}H_t(\al'_\eta)=0.
\eee

The assumptions on $f$ imply that for $\eta$ small enough
$|f''_\eta(x)|\leq \ep_\eta(\eta\bigwedge|x|)$
for a family $\ep_\eta$ of positive numbers going to $0$ as $\eta\to0$.
Then one can checks that the functions $g_\eta$ and $k_\eta$ satisfy
the following, where $\ep'_\eta\to0$ as $\eta\to0$:
$$|g_\eta(x,y)|\leq\ep'_n\phi_2(x)\phi_1(y),\qquad
k_\eta(x,y)^2\leq\ep'_n\phi_2(x)\phi_2(y).$$
We deduce that, for some other family $\ep''_\eta$ still going to $0$
as $\eta\to0$, we have $|\al_\eta|\leq\ep''_n\phi_1$ and
$\al'_\eta\leq\ep''_n\phi_2$. Then (\ref{33}) and (\ref{Conv2}) yield
the first and the second parts of (\ref{95}), respectively. Hence (b)
is completely proved.
\vst

\ni\it Step 4. \rm Finally we prove (c).
When $f=g$ with $g(x)=x^2$, this is Theorem 6.1--(b) of \cite{JP}: the process
$\overline{Z}_n$ of this theorem is indeed $\frac12 Z^n(g')$; this
theorem is proved when $\De_n=1/n$, but the
proof works equally well for any sequence $\De_n\to0$. Moreover a
close look at that proof shows that it also ensures the stable
convergence in law of the pair $(Z^n(g),Z^n(h))$ toward
$(Z(g')+c\sqrt{2}~W',Z(h'))$ (for the Skorokhod topology on the set of
$\R^2$--valued c\`adl\`ag functions), as soon as $h$ is $C^1$ and
vanishes around $0$. Then the sums $Z^n(g)+Z^n(h)$ stably converge
in law to $Z(g')+c\sqrt{2}~W'+Z(h')$. Then if $f\in\ea_2\cap C^1$ the
function $h=f-g$ is as above, hence we get the result. \qed

\subsection{Proof of Theorem \ref{T9}.}

Let us first consider the case where $J'$ and $J''$ are empty, so we
only have components of the first type. We can then reproduce the
previous proof, with a function $f$ which is multidimensional with
components $f_j$. When $f$ (that is, all $f_j$'s) vanishes on a
neighborhood of $0$, we have (\ref{96}) with
$f'=(f'_1,\cdots,f'_d)$ and the variables $\widetilde{R}'^n_p$
depend on $j$, but with the same $R'^n_p$ for all $j$. Then we
deduce from (\ref{93}) the stable convergence in law toward $Z(f')=
(Z(f'_1),\cdots,Z(f'_d))$ with the same sequence $U'_p$. Next, Step 3
of the previous proof is performed component by component, hence there
is nothing to change, and the result is proved.

In a second step we consider the case where $J$ is empty. Since the
limit is then continuous, the conditions for convergence given in Lemmas
\ref{L4} and \ref{LL1} are componentwise, except for the second one in
(\ref{Conv1}). Moreover for each component these conditions are
satisfied, as we have shown in the proof of Theorem \ref{T7}. Hence
the only thing to show is the relevant extension of (\ref{74}). To
write it, we set
$$u_n(j)=\left\{\begin{array}{ll}
\De_n^{1/2-r(j)/2}\quad&\mbox{if }~j\in J'\\[2mm]
\De_n^{1/2}\quad&\mbox{if }~j\in J''\end{array}\right.\qquad
\ze_n^j=\left\{\begin{array}{ll}
\phi_{r(j)}(X_{\De_n})\quad&\mbox{if }~j\in J'\\[2mm]
f_j(X_{\De_n}/\sqrt{\De_n})\quad&\mbox{if }~j\in J''\end{array}\right.$$
Then we need to show that
$$\frac1{\De_n}~u_n(j)u_n(k)\left(\E(\ze_n^j\ze_n^k)-
\E(\ze_n^j)\E(\ze_n^k)\right)~\to~c'^{jk},$$
where $c'$ is given by (\ref{17}). But this is a trivial consequence
of (\ref{33}), (\ref{35}) and (\ref{36'}).

It remains to consider the general case. Let us denote $Y^n(1)
=(Y^{n,j})_{j\in J}$ and $Y^n(2)=Y^{n,j})_{j\in J'\cup J''}$.
These two sequences converge stably in law to
their limits, say $Y(1)$ and $Y(2)$. Moreover $X^{(n)}$ (the
discretized process) has the same jumps than $Y^n(1)$, and converges
(pointwise) to $X$, and all jumps of $Y(1)$ are jumps of $X$, and
$Y(2)$ is continuous: so in fact the pairs $(X^{(n)},Y^n(1))$
and $(X^{(n)},Y^n(2))$ converge stably in law to
$(X,Y(1))$ and $(X,Y(2))$ respectively. Using once more the continuity
of $Y(2)$ we deduce the tightness of the sequence
$Z^n=(X^{(n)},Y^n(1),Y^n(2))$.

Take any subsequence $Z^{n_k}$ which converges in law to a limit
$Z'=(X',Y'(1),Y'(2))$. This limit is necessarily a L\'evy process, and
$(X',Y'(1))$ and $(X',Y'(2))$ have (separately) the same laws as $(X,Y(1))$ and
$(X,Y(2))$: hence $X'$ and $Y'(2)$ are independent, and $Y'(1)$ has no
continuous Gaussian part and $Y'(2)$ is
continuous. By a well known result on multidimensional L\'evy
processes, this implies that $Y'(1)$ and $Y'(2)$ are also independent. In
other words the law of $Z'$ is the product of the laws of $(X,Y(1))$
and of $Y(2)$, and it follows that the original sequence $Z^n$ converges in
law to the process $Z=(X,Y(1),Y(2)$ where $Y(2)$ is independent of
$(X,Y(1))$.

It remains to apply once more the pointwise convergence $X^{(n)}\to
X$, plus the fact that the $\si$--field $\f$ is $\si(X_t:t\geq0)$: we
then deduce from the convergence of $Z^n$ to $Z$, that $Y^n$ converges
stably in law to $Y=(Y(1),Y(2)$ with $Y(2)$ independent of $(X,Y(1)$
and $Y(1)=(Z(f_f))_{j\in J}$, and we are done. \qed

\section{Proofs of the theorems about $\BV^n(f)$ and $\BV'^n(f)$}\label{MT}
\setcounter{equation}{0}
\renewcommand{\theequation}{\thesection.\arabic{equation}}

\nib Proof of Theorem \ref{T11}. \rm (i,ii) It is
obviously  enough to prove the convergence of
$Y^n_t=\frac1{\sqrt{T_n}}(\BV^n(f)_t-H_{\De_n}(f)[nt])$. Now $Y^n_t$ is
the sum of $[nt]$ i.i.d. variables distributed as $\ze_n-\E(\ze_n)$,
with $\ze_n=f(X_{\De_n})/\sqrt{T_n}$. We have $|\ze_n|\leq
C/\sqrt{T_n}\to0$ and $\E(\ze_n^2)=\Ga_{\De_n}(f)/T_n$. Therefore
Lemma \ref{L12}--(a,b) implies that $n\E(\ze_n^2)$ converges to
$F(f^2)$ if $f$ is as described in (i), and to
$F(f^2)+c(\mu_2-\mu_1^2)$ if $f\in\ea_1^b\cap C^{0,F}$: the result
follows from a standard CLT.

(iii) Let $f\in\ea^b_r$ and $r<1$ and either $c>0$ or $2r\notin I$. By
Lemma \ref{L12}--(c) we have $\Ga_{\De_n}(f)/\De_n\to\infty$, hence a
fortiori $u_n:=1/\sqrt{n\Ga_{\De_n}(f)}\to0$, and exactly as above
it is enough to prove the convergence of
$Y^n_t=u_n(\BV^n(f)_t-H_{\De_n}(f)[nt]$, and $Y^n_t$ is
the sum of $[nt]$ i.i.d. variables distributed as $\ze_n-\E(\ze_n)$,
with $\ze_n=u_nf(X_{\De_n})$. We have $|\ze_n|\leq
Cu_n\to0$ and $\E(\ze_n^2)=1/n$, so we conclude as in (i).
\vsc

\nib Proof of Theorem \ref{T13}. \rm We can replace $t$ by
$[nt]/n$ in the centering term of $Y^n$.
Then the proof of the first claim goes exactly as for Theorem \ref{T7}, except
that we use the usual $q$--dimensional CLT and we do not need
something like (\ref{75}). For the last claim, we use (\ref{34}) and
(\ref{36}). \qed
\vsc

\nib Proof of Theorem \ref{T10}. \rm (i) In all cases we deduce
from Theorem \ref{T11}--(i,ii) that
$$\frac1{T_n}\left(\BV^n(f)_t-nH_{\De_n}(f)t\right)=\frac1{T_n}~\!\BV^n(f)_t
-\frac t{\De_n}~\!H_{\De_n}(f)~\toucp~0.$$
Then (\ref{Conv2}) immediately gives the
result.

(ii) and (iii): The same type of arguments, based on Theorem
\ref{T13} and (\ref{33}) and (\ref{35}), gives the result. \qed


\begin{thebibliography}{99}




\bibitem{AJ} A\"it Sahalia Y. and J. Jacod (2005): Volatility
estimators for discretely sampled L\'evy processes.
Preprint.

\bibitem{ABD}
Andersen, T.~G., T.~Bollerslev, and F.~X. Diebold (2005): Parametric
and nonparametric measurement of volatility.  In Y.~A{\"{\i}}t-Sahalia
and L.~P. Hansen (Eds.), {\em Handbook of
  Financial Econometrics}. Amsterdam: North Holland. Forthcoming.

\bibitem{BS}
Barndorff-Nielsen, O.~E. and N.~Shephard (2003): Realised power
variation and stochastic volatility.  {\em Bernoulli\/}~{\bf 9}, 243--265.
Correction published in pages 1109--1111.

\bibitem{BGJPS}
Barndorff-Nielsen, O.~E., S.~E. Graversen, J.~Jacod, M.~Podolskij and
  N.~Shephard (2004): A central limit theorem for realised power and
  bipower variations of continuous semimartingales. To appear.

\bibitem{BSW}
Barndorff-Nielsen, O.~E., N. Shephard, and M. Winkel (2005): Limit
theorems for multipower variation in the presence of jumps. Preprint.


\bibitem{JS}
Jacod, J. and A. Shiryaev (2003): {\it Limit Theorems for Stochastic
Processes}, 2d ed., Springer-Verlag: Berlin.

\bibitem{JP}
Jacod, J. and Protter, P. (1998):
Asymptotic error distributions for
the Euler method for sto\-chas\-tic differential equations. {\it
Ann. Probab.}, {\bf 26}, 267-307.


\bibitem{J}
Jacod, J. (2004):
The Euler scheme for L\'evy driven stochastic differential equations:
limit theorems. {\it Ann. Probab.}, {\bf 32}, 1830--1972.


\bibitem{JJM}
Jacod J., A. Jakubowski, J. M\'emin (2003): On asymptotic error in
discretization of processes. {\it Annals Probab.}, {\bf 31}, 592--608.


\bibitem{L}
L\'epingle, D. (1976): La variation d'ordre $p$ des semimartingales.
{\it Z. f\"ur Wahr. Th.}, {\bf 36}, 285--316.


\bibitem{Ma}
Mancini, C. (2001): Disentangling the jumps of the diffusion in a geometric
jumping Brownian motion. {\it Giornale dell'Instituto Italiano
  degli Attuari} {\bf LXIV} 19--47.

\bibitem{W}
Woerner,J. (2004): Power and multipower variation: inference for high
frequency data. Preprint.



\end{thebibliography}
\end{document}